\documentclass[12pt]{article}

\usepackage{numbysec}
\usepackage{amsmath}
\usepackage{amsthm}
\usepackage{amsfonts}
\usepackage{amssymb}
\usepackage{latexsym}
\usepackage{float}
\usepackage{upref}
\restylefloat{figure}
\usepackage{graphicx}
\usepackage{color}

\hsize \textwidth
\advance \hsize by -\marginparwidth
\evensidemargin \oddsidemargin

\def\R{\mathbb{R}}
\def\T{\mathbb{T}}
\newcommand{\RNum}[1]{\text{\uppercase\expandafter{\romannumeral #1\relax}}}

\newtheorem{thm}{Theorem}[section]
\newtheorem{lem}[thm]{Lemma}
\newtheorem{cor}[thm]{Corollary}
\newtheorem{prop}[thm]{Proposition}

\theoremstyle{assumption}
\theoremstyle{definition}
\newcommand{\eps}{\varepsilon}
\newcommand{\la}{\lambda}

\newcommand{\farc}{\frac}
\newcommand{\Rm}{{\mathbb R}}

\numberbysection
\begin{document}

\title{
Global regularity for the fractional Euler alignment system}

\author{Tam Do\thanks{Department of Mathematics, Rice University, Houston, TX 77005, USA.
E-mail: Tam.Do@rice.edu.}
\and Alexander Kiselev\thanks{Department of Mathematics, Rice University, Houston, TX 77005, USA.
E-mail: kiselev@rice.edu.}
\and
Lenya Ryzhik\thanks{Department of Mathematics, Stanford University, Stanford, CA 94305, USA. Email: ryzhik@stanford.edu.}
\and
Changhui Tan\thanks{Department of Mathematics, Rice University, Houston, TX 77005, USA.
E-mail: ctan@rice.edu.}
}

\maketitle

%%%%%%%%%%%%%%%%%%%%%%%%%%%%%%%%%%%%%%%
% Abstract
%%%%%%%%%%%%%%%%%%%%%%%%%%%%%%%%%%%%%%%
\begin{abstract}
%Nonlinear partial differential equations modeling collective behaviors have been the focus
%of much recent research.
We study a pressureless Euler  system with a non-linear
density-dependent alignment term, originating in the Cucker-Smale
swarming models. The alignment term is dissipative in the sense that it
tends to equilibrate the velocities.  Its density dependence  is natural:
the alignment rate  increases in the areas of high density due to
species discomfort.
The diffusive term has the order of a fractional Laplacian~$(-\partial_{xx})^{\alpha/2}$,~$\alpha\in(0,1)$.
The
corresponding Burgers equation with a linear dissipation of this type
develops shocks in a finite time. We show that
the alignment nonlinearity enhances the dissipation,
and the solutions are globally regular for all $\alpha\in(0,1)$.
To the best of our knowledge, this is the first example of  such
regularization due to the non-local nonlinear modulation of dissipation.
\end{abstract}

%%%%%%%%%%%%%%%%%%%%%%%%%%%%%%%%%%%%%%%
% Title and Table of contents
%%%%%%%%%%%%%%%%%%%%%%%%%%%%%%%%%%%%%%%
%\maketitle
%\vspace*{-0.8cm}
%\tableofcontents

\section{Introduction}

\subsubsection*{The Cucker-Smale model}

Modeling of the self-organized collective behavior, or swarming, has
attracted a large amount of attention over the last few years. Even an attempt at a brief
review of this field is well beyond the scope of this introduction, and we refer to the recent
reviews~\cite{carrillochoiperez,choihaperez,Vicsek2012}.
A remarkable phenomenon commonly observed in biological
systems is  flocking, or velocity alignment by near-by individuals.
One  of the early
flocking models, discrete in time and two-dimensional,  is commonly referred to as the Vicsek model: the angle $\theta_i(t)$
of the velocity of $i$-th particle satisfies
\begin{equation}\label{jan502}
\theta_i(t+1)=\frac{1}{|{\cal N}_i(t)|}\sum_{j\in{\cal N}_i(t)}\theta_j(t)+
\eta \Delta\theta.
\end{equation}
Here, ${\cal N}_i(t)=\{j:~|x_i(t)-x_j(t)|\le r\}$, with some $r>0$ fixed,
$\Delta\theta$ is a uniformly distributed random variable in $[-1,1]$,
and $\eta>0$ is a parameter measuring the strength of the noise.
This model  preserves the modulus of the particle velocity
and only affects its direction. First via numerical simulations and then by mathematical tools, it has been shown
that this model has a rich behavior, ranging from flocking when~$\eta$ is small,
to a completely chaotic motion for large~$\eta$, with a phase transition at
a certain critical value $\eta_c$.

A natural generalization of the Vicsek model was introduced by F. Cucker
and S. Smale~\cite{cucker2007emergent}:
\begin{equation}\label{jan602}
\dot{x}_i=v_i,\quad \dot{v}_i=\frac{1}{N}\sum_{j=1}^N\phi(|x_i-x_j|)(v_j-v_i).
\end{equation}
Here, $\{x_i, v_i\}_{i=1}^N$ represent,
respectively, the locations and the velocities of the
agents. Individuals align their velocity to the neighbors,
with the interaction strength characterized by a non-negative influence
function~$\phi(x)\ge 0$.
The relative influence is typically taken as a decreasing function of the distance
between individuals. An important flexibility of the Cucker-Smale model is that it
both does not
impose the constraint on the velocity magnitude and allows to analyze the behavior based
on the decay properties of the kernel $\phi(r)$. One of the main results of the Cucker-Smale
paper was that, roughly,
provided that $\phi(r)$ decays slower than~$r^{-1}$ as~$r\to+\infty$, then
all velocities $v_i(t)$ converge to a common limit $\bar v(t)$,
and the relative particle positions
$x_i(t)-x_j(t)\to \bar x_{ij}$ also have a common limit -- the particles form a
swarm moving with a uniform velocity. This is what we would call a global flocking:
all particles move with nearly identical velocities.

One potential shortcoming of the Cucker-Smale model is that an "isolated clump"
of particles may be more affected by "far away" large mass than by its
own neighbors. Essentially, the dynamics inside a small clump would
be suppressed by the presence of a large group of particles "far away".
This can be balanced by a different kind of averaging,
rather than simple division by~$N$ in (\ref{jan602}),
as was done by S. Motsch and E. Tadmor in~\cite{motsch-tadmor11}:
\begin{equation}\label{jan604}
\dot{x}_i=v_i,\quad \dot{v}_i=\frac{\lambda}{\Phi_i}
\sum_{j=1}^N\phi(|x_i-x_j|)(v_j-v_i),~~\Phi_i=\sum_{k=1}^N\phi(|x_i-x_k|),
\end{equation}
with some $\lambda>0$.
This modification reinforces the local alignment over the long distance
interactions.

\subsubsection*{A kinetic Cucker-Smale model}

Kinetic models are also commonly used to describe the collective behavior when the number of particles
is large, in terms
of the particle density $f(x,v,t)$, with $x\in\Rm^d$, $v\in\Rm^d$.
A kinetic limit of the Cucker-Smale model was obtained by S.-Y. Ha and E. Tadmor in~\cite{ha2008particle},
%When the number of particles tends to infinity, it is convenient to use
%a macroscopic description
%Let us recall the Liouville equation
%for a single particle obeying the dynamics
%\begin{equation}\label{jan504}
%\dot x=v,~~\dot v=G(t,x,v).
%\end{equation}
%Then the particle density obeys the first order PDE, known as the Liouville equation:
as a nonlinear and non-local kinetic equation
\begin{equation}\label{jan506}
f_t+v\cdot\nabla_x f+\nabla_v\cdot(L[f] f)=0,
\end{equation}
with
\begin{equation}\label{jan508}
L[f](x,v,t)=\int_{\Rm^{2d}}\phi(x-y)(v'-v)f(y,v',t)dv'dy.
\end{equation}
Together, (\ref{jan506})-(\ref{jan508}) give a nonlinear kinetic version of the Cucker-Smale system.
%Ha and Tadmor have shown that it admits classical solutions for the kernels originally considered by
%Cucker and Smale:
%\begin{equation}\label{jan510}
%\phi(x)=\farc{C}{(1+|x|^{2})^\beta}.
%\end{equation}
% It was shown in~\cite{ha2008particle} that its solutions exhibit global flocking, in the sense that
% \begin{equation}\label{jan512}
% \int_{\Rm^{2d}}|v-u_c|^2f(t,x,v)dxdv\to 0~~\hbox{as $t\to+\infty$},
% \end{equation}
% under the assumption that $\phi(r)$ decays sufficiently slowly as $r\to+\infty$.
% Here, $u_c$ is the mean velocity:
% \begin{equation}\label{jan514}
% u_c=\frac{1}{M_0}\int vf(t,x,v)dxdv,~~~ M_0=\int_{\Rm^{2d}}f(t,x,v)dxdv.
% \end{equation}
% Both the total mass $M_0$ and the mean velocity $u_c$ are preserved in time.

It was shown in \cite{carrillo2010asymptotic} that its solutions
exhibit global flocking, in the sense that the size of the support in $x$
\[
S(t)=\sup\{|x-y|~:~(x,v),(y,v')\in\text{supp}(f(\cdot,\cdot,t))\}
\]
remains uniformly bounded in time, and
the support in $v$ shrinks:
\begin{equation}\label{jan1201}
V(t)=\sup\{|v-v'|~:~(x,v),(y,v')\in\text{supp}(f(\cdot,\cdot,t))\}\to
0\hbox{ as $t\to+\infty$,}
\end{equation}
under the assumption that $\phi(r)$ decays slower than $r^{-1}$ as
$r\to+\infty$.
A similar result was obtained in \cite{tan2017discontinuous} for the kinetic
Motsch-Tadmor system.

A kinetic model that combines the features of the Cucker-Smale and Motsch-Tadmor
models was proposed in a paper by T. Karper, A. Mellet and K. Trivisa~\cite{karper-mellet-trivisa}:
\begin{equation}\label{jan608}
f_t+v\cdot\nabla_x f+\nabla_v\cdot(L[f] f)+\lambda\nabla_v\cdot((u(x,t)-v)f)=\Delta_v f,
\end{equation}
with $L[f]$ as in (\ref{jan508}), $\lambda>0$, and
the local average velocity $u(t,x)$ defined as
\begin{equation}\label{jan610}
u(x,t)=\frac{1}{\rho(t,x)}\int_{\Rm^d}v f(x,v,t)dv,~~
\rho(x,t)=\int_{\Rm^d} f(x,v,t)dv.
\end{equation}
The Laplacian in the right side of (\ref{jan608}) takes into account the
possible Brownian noise in the velocity.

One should also mention a large body of literature on the kinetic versions of the Vicsek model and its modifications,
and their hydrodynamic limits:
see~\cite{carlen-carvalho-degond-wennberg,degond-dimarco-mac-wang,degond-frouvelle-liu,degond-frouvelle-raoul,degond-liu-mosch-panferov}
and references therein.

\subsubsection*{An Euler alignment model}

The kinetic Cucker-Smale
model can be further "macroscopized" as a hydrodynamic model for the local density $\rho(t,x)$ and local average velocity $u(t,x)$ defined in
(\ref{jan610}).
%\begin{equation}\label{jan516}
%\rho(t,x)=\int_{\Rm^d}f(t,x,v)dv,~~u(t,x)=\int_{\Rm^d} v f(t,x,v)dv.
%\end{equation}
The standard formal derivation of the hydrodynamic limit for nonlinear kinetic equations
often relies on a
(often hard to justify) moment closure
procedure.
An alternative is to consider the "monokinetic" solutions of (\ref{jan506})-(\ref{jan508})
of the form
\begin{equation}\label{jan518}
f(x,v,t)=\rho(t,x)\delta(v-u(x,t)).
\end{equation}
In a sense, this is a "local alignment" (as opposed to global flocking) ansatz -- the particles move
locally with just a single velocity but the velocity does vary in space.
Inserting this expression into (\ref{jan506})-(\ref{jan508})
gives the Euler alignment system, which we write in one dimension as
\begin{eqnarray}\label{oct1002}
&&  \partial_t\rho+\partial_x(\rho u)=0, \\
&& \partial_t(\rho u)+\partial_x(\rho u^2)=
\int_\R\phi(x-y){(u(y,t)-u(t,x))} \rho(y,t)\rho(x,t)dy.\label{oct1004}
\end{eqnarray}
The presence of the density $\rho$ under the integral in the right
side of \eqref{oct1004} has a
very reasonable biological interpretation: the alignment effect between the individual agents becomes stronger
where the density is high (assuming that the interaction kernel $\phi$
is localized).   As far as a rigorous derivation of
the hydrodynamic limit is concerned, the aforementioned
paper~\cite{karper-mellet-trivisa} derives the hydrodynamic limit
starting from the "combined" Cucker-Smale-Motsch-Tadmor
kinetic system (\ref{jan608}):
\begin{eqnarray}\label{jan618}
&&  \partial_t\rho+\partial_x(\rho u)=0, \\
&& \partial_t(\rho u)+\partial_x(\rho u^2)+\partial_x\rho=
\int_\R\phi(x-y){(u(y,t)-u(x,t))} \rho(y,t)\rho(x,t)dy.\label{jan620}
\end{eqnarray}
This system has an extra term $\partial_x\rho$ in
the left side of (\ref{oct1004}) that can be thought of as pressure,
with the constitutive law $p(\rho)=\rho$. The pressure appears as a result of the
balance between the   local interaction term in the left side of (\ref{jan608}) and the Laplacian
in the right side.
In particular, the starting point of the derivation is not the single
local velocity ansatz (\ref{jan518})
but its smooth Maxwellian version (setting~$\lambda=1$ in~(\ref{jan608}) for convenience)
\begin{equation}\label{jan612}
f(x,v,t)=\rho(x,t)\exp\Big(-\frac{(v-u(x,t))^2}{2}\Big),
\end{equation}
together with the assumption that the interaction is weak: $\phi\to\eps\phi$,
and a large time-space rescaling~$(t,x)\to t/\eps,x/\eps$.

Another version of the Euler equations as a model
for swarming has been proposed in~\cite{levine-rappel-cohen},
and formally justified in~\cite{bertozzi-chayes}:
\begin{eqnarray}\label{jan614}
&&\partial_t\rho+\partial_x(\rho u)=0,\\
&&\partial_t(\rho u)+\partial_x(\rho u^2)+\partial_x\rho=\alpha\rho u-
\beta \rho|u|^2u-
\int_\R\nabla V(x-y)\rho(y,t)\rho(x,t)dy.\nonumber
\end{eqnarray}
The key difference between models like \eqref{jan614} and the ones we consider here
is the absence of the regularizing term $u(t,y)-u(t,x)$
in the right side, so one does not expect the regularizing effect of the
interactions that we will observe here.

\subsubsection*{The Euler alignment system for Lipschitz interaction kernels}

When particles do not interact, that is, $\phi(x)\equiv 0$,
the system (\ref{oct1002})-(\ref{oct1004}) is simply the pressure-less Euler equations. In particular,
in that case,~(\ref{oct1004}) is   the Burgers equation:
\begin{equation}\label{jan520}
u_t+uu_x=0.
\end{equation}
Its solutions
develop a shock singularity in a finite time if the initial
condition $u_0(x)$ has a point where $\partial_xu_0(x)<0$. In particular,
if $u_0(x)$ is periodic and not identically equal to a constant,
then $u(x,t)$
becomes discontinuous in a finite time.
The function~$z(x,t)=-u_x(x,t)$   satisfies the continuity
equation
\begin{equation}\label{jan802}
z_t+(zu)_x=0,
\end{equation}
and becomes infinite at the shock location.

The singularity
in the Burgers equation does not mean
that there is a singularity in the solution of kinetic equation: it only means that
the ansatz~(\ref{jan518}) breaks down,
and we can not associate a single velocity to a given position. This is a version
of "a shock implies no local alignment". To illustrate this point, consider the solution of
free transport equation
\begin{equation}\label{jan804}
f_t+vf_x=0,
\end{equation}
with the initial condition $f_0(x)=\delta(v+x)$. The solution of the kinetic
equation is
\begin{equation}\label{jan806}
f(x,v,t)=f_0(x-vt,v)=\delta(v+x-vt),
\end{equation}
hence the ansatz (\ref{jan518}) fails at $t=1$. This is the time when the corresponding Euler equation
\begin{equation}\label{jan902}
u_t+uu_x=0,
\end{equation}
with the initial condition $u(0,x)=-x$, develops a shock: $u(x,t)=-x/(1-t)$.

The integral
term in the right side of (\ref{oct1004}) has a dissipative nature
when $\phi\not\equiv 0$: it tries to
regularize the velocity discontinuity.
When the function $\phi(x)$ is Lipschitz, this system has been investigated
in~\cite{carrillo2016critical} and~\cite{tadmor2014critical} that
show two results. First, a version of global flocking: if $\phi$ decays slower than $|x|^{-1}$
at infinity, and the solution remains smooth for all
$t\ge 0$ and the initial density $\rho_0$ is compactly
supported, then the support $S_t$ of~$\rho(t,\cdot)$ remains uniformly bounded in
time, and
\begin{equation}\label{jan522}
\sup_{x,y\in S_t}|u(x,t)-u(y,t)|\to 0\hbox{ as $t\to+\infty$.}
\end{equation}
An improvement in global regularity compared to the Burgers
equation (\ref{jan520}) was also obtained in~\cite{carrillo2016critical} and~\cite{tadmor2014critical}.
As we have mentioned,
solutions of the latter become discontinuous in a finite time provided there is a
point $x\in\Rm$ where the initial condition $u_0(x)$ has a negative derivative:
$\partial_x u_0(x)<0$. On the other hand, solutions of the Euler alignment
equations remain
regular  for initial data such that
\begin{equation}\label{oct1408}
\partial_xu_0(x)\ge -(\phi\star \rho_0)(x)\hbox{ for all $x\in\Rm$,}
\end{equation}
while
the solution blows up in a finite time if there exists~$x_0\in\Rm$ such that
\begin{equation}\label{oct1410}
\partial_xu_0(x)< -(\phi\star \rho_0)(x).
\end{equation}
Thus, the presence of the dissipative term in~(\ref{oct1004}) leads to global regularity for some initial data
that blows up for the Burgers equation: the right side of (\ref{oct1408})
may be negative.
%Note that the presence of the density $\rho$ under the integral on the right hand side of \eqref{oct1004} has a
%very reasonable biological interpretation: the alignment effect between the individual agents becomes stronger
%where the density is high.
%Thus, the transition between the initial conditions that lead
%to regular solutions and those that lead to shocks is shifted by the dissipative term in~(\ref{oct1004}) in favor of the
%global regularity.
However, a Lipschitz interaction kernel $\phi(x)$
arrests the shock singularity for the Euler alignment
equations only for some initial conditions.

\subsubsection*{Singular alignment kernels}

Our interest is in singular interaction kernels of the form~$\phi(x)=|x|^{-\beta}$, with $\beta>0$.
One reason to consider such kernels is to strengthen the effect of the local
interactions compared to the effect of "far-away" particles, in the spirit
of the Motsch-Tadmor correction.
The well-posedness of the finite number of particles Cucker-Smale system with such interactions is a delicate
issue -- the difficulty is in either ruling out the possibility of particle collisions, or understanding the
behavior of the system at and after a collision.
This problem was addressed in~\cite{peszek2014,peszek2015}
for~$\beta\in(0,1)$ -- it was shown that in this range, particles may get stuck together
but a weak solution of the ODE system can still be defined. When $\beta\ge 1$, a set of initial conditions that has no particle
collisions was described in~\cite{ahn-choi-lee2012}. The absence of collisions was proved very recently for general initial configurations
in~\cite{carrillo-choi-mucha-peszek2016}.
As far as flocking is concerned, unconditional flocking was proved in~\cite{ha-liu2009} for $\beta\in(0,1)$,
while for $\beta\ge 1$ there are initial configurations that do not lead to global flocking -- the long distance
interaction is too weak.
The well-posedness of the kinetic Cucker-Smale system for~$\beta \in (0,1/2)$ was established
in~\cite{mucha-peszek2015}.

We consider here the alignment kernels $\phi(x)$ with $\beta>1$:
\begin{equation}\label{oct1702}
\phi_\alpha(x)=\farc{c_\alpha}{|x|^{1+\alpha}},
\end{equation}
with $\alpha>0$.
In particular, the decay of $\phi(x)$
at large $|x|$ is faster than the $1/|x|$ decay required for the Cucker-Smale and other
proofs of flocking. It is compensated by a very strong alignment for $|x|\to 0$.
The constant $c_\alpha$ is chosen so that
\[
\Lambda^\alpha f=c_\alpha\int_\Rm
\farc{f(x)-f(y)}{|x-y|^{1+\alpha}}dy,~~\Lambda=(-\partial_{xx})^{1/2}.
\]
Then the strong form of the
Euler alignment  system is
\begin{eqnarray}
&&\partial_t\rho+\partial_x(\rho u)=0\label{eq:mainrho}\\
&&\partial_tu+u\partial_xu=c_\alpha\int_\R\frac{u(y,t)-u(x,t)}{|y-x|^{1+\alpha}}
\rho(y,t)dy.\label{eq:mainu}
\end{eqnarray}

% We note that another interesting
% case is a ``very localized" $\phi(x)$ such
% as $\phi_\eps(x)=\eps^{-3}\phi(x/\eps)$, with~$\eps\ll 1$ and a smooth
% compactly supported function $\phi(x)$. The hyperbolic scaling
% $(x,t)\to(x/\eps, t/\eps)$ together with the formal limit
% $\eps\to 0$ leads the compressible Navier-Stokes equations
% \begin{eqnarray}\label{oct1302}
% &&  \partial_t\rho+\partial_x(\rho u)=0, \\
%  && \partial_t(\rho u)+\partial_x(\rho u^2)= \nu(\mu(\rho)u_x)_x,\label{oct1304}
% \end{eqnarray}
% with the viscosity $\mu(\rho)=\rho^2$, and
% \[
% \nu=\farc 12\int_\Rm\phi(y)y^2dy.
% \]
% The compressible Navier-Stokes equations with a density-dependent viscosity
% have been extensively studied in the past, with the main focus
% on the issue of a vacuum region, where~$\rho$ vanishes.
% We are concerned here with a different problem of
% preventing a shock, and
% we should mention the paper by A.~Mellet and A. Vasseur~\cite{melletvasseur}
% that studies the global
% existence of regular solutions to (\ref{oct1302})-(\ref{oct1304})
% with isentropic pressure and~$\mu(\rho)$ that can vanish at $\rho=0$
% as $\mu(\rho)\sim\rho^m$, $0\le m\le 1/3$.

Let us first compare the
Euler alignment  system~(\ref{eq:mainrho})-(\ref{eq:mainu})
to the  Burgers equation with a fractional dissipation
\begin{equation}\label{eq:fractalBurgers}
\partial_tu+u\partial_xu=-\Lambda^\alpha u,
\end{equation}
obtained by formally setting $\rho(t,x)\equiv 1$ in (\ref{eq:mainu}) and dropping
(\ref{eq:mainrho}) altogether. This neglects the nonlinear mechanism of the
dissipation.
Global regularity of the solutions of the fractional Burgers equation
has been studied in~\cite{kiselev2008blow}. One can distinguish three
regimes:
first, when $\alpha>1$, the dissipative term in the right side has a higher order
derivative than the  nonlinear term in the left side. This is the sub-critical
regime: the dissipation
dominates the nonlinearity, and global existence of the strong solutions can
be shown in a reasonably straightforward manner
using the energy methods. On the other hand, when $0<\alpha<1$, the
dissipation is too weak to compete with the nonlinear term, which has a higher
derivative,
and solutions with smooth initial conditions may develop a shock, as
in the inviscid case. The critical case is $\alpha=1$ when the dissipation and the
nonlinearity
contain derivatives of the same order. One may expect that then the nonlinearity
may win
over the dissipation for some
large data. This, however, is not the case: solutions with smooth initial
conditions remain regular globally in time.
The proof of the global regularity when~$\alpha=1$ is much less straightforward
than for~$\alpha>1$
and does not rely solely on the energy methods.

One may hope that the nonlinearity in the dissipative term in the right side of (\ref{eq:mainu})
is actually beneficial, compared to the fractional Burgers equation
(\ref{eq:fractalBurgers}).
Indeed, on the qualitative level, as the shock would form, the density~$\rho$ would
be expected to increase near the
point of the shock. This, in turn, would increase
the dissipation in (\ref{eq:mainu}), moving the problem from "like a super-critical
Burgers" to "like a sub-critical Burgers".
This intuition,  however,
may be slightly misleading: for instance, as we will see, strengthening the
dissipation by
increasing~$\alpha$ does not appear to make the problem any easier, or change its critical character. The
competition between the Burgers nonlinearity  in the left
side of~(\ref{eq:mainu}) and the nonlinear dissipation in the right side
is rather delicate.

The aforementioned results of~\cite{carrillo2016critical,tadmor2014critical} may lead  to
a conjecture that a dissipation term involving the
convolution kernel $\phi \notin L^1$, as in \eqref{eq:mainu},
should lead to global regularity.  However,  this is far from obvious.
The global regularity argument of~\cite{carrillo2016critical,tadmor2014critical}  uses two ingredients:  first, if initially
\begin{equation}\label{jan624}
\partial_x u_0 +
\phi \star \rho_0 \geq 0,
\end{equation}
for all $x\in\Rm$ then
\begin{equation}\label{jan626}
\partial_x u + \phi \star \rho \geq 0
\end{equation}
for all $x\in\Rm$ and $t\ge 0$. Second,  an $L^\infty$-bound
on $\rho$ is established. When $\phi$ is
an~$L^1$-function, one deduces
a lower bound $\partial_x u > -C_0$,
which is crucial  for global regularity. One may combine
an argument of~\cite{carrillo2016critical} with the Constantin-Vicol
nonlinear maximum principle to establish the $L^\infty$-bound for $\rho$ in
our case, as well.
However, in our
case, the analogous inequality to (\ref{jan626})~is
\begin{equation}\label{jan628}
\partial_x u - \Lambda^\alpha \rho \geq 0.
\end{equation}
This fails to give the required
lower bound on $\partial_x u$
based on just the~$L^\infty$ control of $\rho$, and the global regularity
does not follow easily from the uniform bound on the density.
Instead, we have to deploy
a much subtler argument
involving both upper and lower bounds on the density and
a non-trivial modification of the modulus of continuity technique of
\cite{kiselev2007global}.

\subsubsection*{The main result}

We consider here the Euler alignment system (\ref{eq:mainrho})-(\ref{eq:mainu})  on the torus $\T$,
for $\alpha\in(0,1)$. In particular, this range of $\alpha$ corresponds to the
supercritical case for the fractional Burgers equation~(\ref{eq:fractalBurgers}). We prove that the
nonlinear, density modulated dissipation qualitatively changes the behavior of
the solutions: instead of blowing up in a finite time, solutions are
globally regular.
%We break the results into three major components.
\begin{thm} \label{thm:global}
For $\alpha\in(0,1)$, the Euler alignment system (\ref{eq:mainrho})-(\ref{eq:mainu})
with periodic smooth initial data $(\rho_0, u_0)$ such that $\rho_0(x)>0$ for all $x\in\T$,
has a unique global smooth solution.
\end{thm}
The regularizing effect of a non-linear diffusion has been observed before, for instance,
in the chemotaxis problems with a nonlinear diffusion --
see \cite{bedrossian2015,bedrossian-kim2013,bedrsossian-rodriguez2014,bedrossian-rodr-berotzzi2011,bertozzi-slepcev2010}.
The main novelties here are that the nonlinearity is non-local, and that, as we will see, increasing $\alpha$ does not,
contrary to a naive intuition, and unlike what happens in the fractional Burgers equation,
strengthen the regularization effect.

%To the best of our knowledge, this effect is completely new and has not been studied before.
%It is also interesting that it extends over the whole range of dissipation parameter $\alpha.$

To explain the ideas behind the result and its proof,
it is convenient to reformulate the Euler alignment system (\ref{eq:mainrho})-(\ref{eq:mainu})  as the following
system for $\rho$ and $G=\partial_xu-\Lambda^\alpha\rho$:
\begin{eqnarray}
&&\partial_t\rho+\partial_x(\rho u)=0,\label{eq:rhobis}\\
&&\partial_tG+\partial_x(Gu)=0,\label{eq:Gbis}
\end{eqnarray}
with the velocity $u$ related to $\rho$ and $G$ via
\begin{equation}
\partial_xu=\Lambda^\alpha\rho+G.\label{eq:uxbis}
\end{equation}
We show in Section~\ref{sec:2} that (\ref{eq:mainrho})-(\ref{eq:mainu})  and (\ref{eq:rhobis})-(\ref{eq:uxbis}) are, indeed, equivalent for regular solutions.
Note that (\ref{eq:uxbis}) only defines $u$ up to its mean, which is determined from the conservation of the momentum:
\begin{equation}\label{jan630}
\int_\T \rho(x,t)u(x,t)dx=\int_\T \rho_0(x)u_0(x)dx.
\end{equation}
%We explain in Section~\ref{sec:2} how the mean is determined in terms of~$\rho(x,t)$ and~$G(x,t)$.
Somewhat paradoxically, (\ref{eq:uxbis}) seems to indicate that increasing the dissipation $\alpha$
makes the velocity more singular in terms of the density rather than more regular.
%The motivation for the introduction of the function~$G$
%can be seen from the global regularity
%conditions~(\ref{oct1408})-(\ref{oct1410}) for the Euler alignment system with a Lipschitz
%influence function.
%

The solutions of (\ref{eq:rhobis})-(\ref{eq:Gbis})  with the initial conditions~$\rho_0(x)$,~$u_0(x)$ such that
\begin{equation}\label{oct1412}
G_0(x)=\partial_x u_0(x)-\Lambda^\alpha\rho_0(x)\equiv 0,
\end{equation}
preserve the constraint $G=0$ for all $t>0$,
and (\ref{eq:rhobis})-(\ref{eq:Gbis})  then reduces to a single equation
\begin{equation}\label{eq:special}
\partial_t\rho+\partial_x(\rho
u)=0,\quad \partial_xu=\Lambda^\alpha\rho,
\end{equation}
that is simpler to analyze.
%
%As seen from the global regularity condition (\ref{oct1408}) for Lipschitz influence kernels, condition (\ref{oct1412}) is natural:
%such solutions are globally regular but are exactly at the boundary for regularity when~$\phi(x)$ is Lipschitz.  Of course,
%(\ref{oct1408}) does not, strictly speaking, make sense as in the present case $\phi(x)=c_\alpha|x|^{-(1+\alpha)}$ is not an
%$L^1$-function, and $\rho_0$ is strictly positive.
Note that (\ref{eq:special}) defines $u(x,t)$ only up to its spatial average --
we assume that it has mean-zero for all $t>0$.
The model (\ref{eq:special}) is interesting in its own right.
When~$\alpha=1$, so that the velocity is  the Hilbert transform of the
density, it was introduced as a 1D vortex sheet model in~\cite{bakermorlet}, and has been extensively studied in
\cite{chae2005finite}  as a 1D model of the~2D
quasi-geostrophic equation. In particular, the global existence of the solution
if~$\rho_0>0$ is proved in~\cite{chae2005finite} using the algebraic properties of the Hilbert transform.
Our results in this paper can
be directly applied to \eqref{eq:special},
and show the global regularity of the solutions for all~$\alpha\in(0,1)$.  The strategy of the regularity
proof here is very different from that
in~\cite{chae2005finite}.
A quintessential feature of (\ref{eq:special}) is that increasing $\alpha$ does not help the dissipation
in its competition with the Burgers nonlinearity. Indeed, the toy model (\ref{eq:special}) can be written as
\begin{equation}\label{oct1416}
\partial_t\rho+(\partial_x^{-1}\Lambda^\alpha\rho)\partial_x\rho=-\rho\Lambda^\alpha\rho.
\end{equation}
Thus, the scalings of the dissipation in the right side and of the nonlinear transport term in the left side are exactly the same,
both in $\rho$ and in $x$, no matter what $\alpha\in(0,1)$ is.
While the proof of global regularity for \eqref{oct1416} is inspired by the nonlocal maximum principle arguments 
of \cite{kiselev2007global,kiselev2008blow}, the nonlinear nature of dissipative term 
necessitates significant changes and new estimates. The upgrade of the proof from global regularity  
of the model equation to the full system is also highly non-trivial and requires new ideas. 

We note that our results can be applied to the case   $\alpha\in(1,2)$, where the
global behavior is the same as for the fractional Burgers equation.
One can also extend our results to influence kernels of the form
\begin{equation}\label{oct303}
\phi(x)=\farc{\chi(|x|)}{|x|},
\end{equation}
with a non-negative smooth compactly supported function $\chi(r)$. This is the analog of
the kernels in~(\ref{oct1702}) for $\alpha=0$.
%Essentially, the only modification required in the proof is
%in the use of the Constantin-Vicol non-local maximum
%principle~\cite{constantin2012nonlinear}.
%{\bf (CT: Not sure whether the energy estimate and modulus
%  argument can be applied easily.)}
%An interested reader may verify that
%the extension is straightforward.
We expect that as soon as the influence kernel is not integrable,
%so that
%condition (\ref{oct1408}) holds automatically,
solutions remain regular.
The proofs of these extensions require some nontrivial adjustments and further technicalities
compared to the arguments in this paper, and will be presented elsewhere.

Our results also lead to global flocking behavior for
(\ref{eq:mainrho})-(\ref{eq:mainu}). 
The periodized influence function
\[
\phi_p(x)=\sum_{m\in\R\backslash\T}\phi(x+m)
\]
has a positive lower bound for all $x\in\T$.
Since the solution is smooth, one can use the argument in
\cite{tadmor2014critical} to obtain asymptotic flocking behavior in the sense
that
\begin{equation}\label{jan1501}
\sup_{x,y\in\T}|u(x,t)-u(y,t)|\to 0\hbox{ as $t\to+\infty$.}
\end{equation}

This paper is organized as follows. In Section~\ref{sec:2} we prove
an a priori $L^\infty$-bound on~$\rho$, that is the key estimate for the regularity
of the solutions, as well as lower bound on $\rho$.
The local well-posedness of the solutions is proved in Section~\ref{sec:3}.
Section~\ref{sec:4}
contains the proof of our main result, Theorem~\ref{thm:global}.
Apppendix~\ref{sec:app} contains the proof of an auxiliary technical  estimate.
Throughout the paper
we denote by $C$, $C'$, etc. various universal constants, and by $C_0$, $C_0'$ etc.
constants that depend only on the initial conditions.

{\bf Acknowledgment.} This work was was partially supported by the NSF grants DMS-1412023 and
DMS-1311903.

\section{Bounds on the density }\label{sec:2}

In this section, we prove the upper and lower bounds on the density $\rho(t,x)$.
The upper bound is uniform in time, and is crucial for the global
regularity. The lower bound will deteriorate in time but will
be sufficient for our purposes.

\subsection{The reformulation of the Euler alignment system}

We first explain how  the Euler alignment system (\ref{eq:mainrho})-(\ref{eq:mainu})
is  reformulated
as ~(\ref{eq:rhobis})-(\ref{eq:Gbis}), as we will mostly use the latter.
We only need to obtain (\ref{eq:Gbis}) for $G$ defined in (\ref{eq:uxbis}).
The idea comes from \cite{carrillo2016critical}.
We apply the operator $\Lambda^\alpha$ to \eqref{eq:mainrho}, and use the identity
\[
u(y)\rho(y)-u(x)\rho(x)=[u(y)-u(x)]\rho(y)+u(x)[\rho(y)-\rho(x)],
\]
to obtain
\begin{eqnarray}\label{oct1704}
&&\partial_t\Lambda^\alpha\rho=
-\partial_x\Lambda^\alpha(\rho u)=%c_\alpha\partial_x\int_\R\frac{u(y)\rho(y)-u(x)\rho(x)}{|y-x|^{1+\alpha}}\rho(y)dy=
c_\alpha\partial_x\int_\R\frac{u(y)-u(x)}{|y-x|^{1+\alpha}}\rho(y)dy
-\partial_x\left(u(x)\Lambda^\alpha\rho\right).
\end{eqnarray}
On the other hand, applying $\partial_x$ to \eqref{eq:mainu}, we get
\begin{equation}\label{oct1706}
\partial_t(\partial_xu)+ \partial_x(u\partial_xu) =
c_\alpha\partial_x\int_\R\frac{u(y)-u(x)}{|y-x|^{1+\alpha}}\rho(y)dy.
\end{equation}
Subtracting (\ref{oct1704}) from (\ref{oct1706}) gives an equation  for the function $G=\partial_xu-\Lambda^\alpha\rho$:
\[
\partial_tG+\partial_x(Gu)=0,
\]
which is (\ref{eq:Gbis}).
%Therefore, system \eqref{eq:main} can be expressed as the coupled
%dynamics of $(\rho, G)$ as follows.
%\begin{subequations}\label{eq:rhoG}
%\begin{align}
%\partial_t\rho+\partial_x(\rho u)=0,\label{eq:rho}\\
%\partial_tG+\partial_x(Gu)=0,\label{eq:G}\\
%\partial_xu=\Lambda^\alpha\rho+G.\label{eq:ux}
%\end{align}
%\end{subequations}
%

Let us comment on how to recover $u$ from \eqref{eq:uxbis}.
Let us denote by
\begin{equation}\label{oct1722}
\kappa=\frac{1}{|\T|}\int_\T\rho(x,t) dx
\end{equation}
the average of $\rho$ in $\T$, which is preserved in time by (\ref{eq:rhobis}), at least as long as $\rho$
remains smooth. Note that $G(x,t)$ has mean zero automatically:
\begin{equation}\label{oct1712}
\int_\T G(x,t)dx=\int_\T G_0(x)dx=0.
\end{equation}
We also define
\begin{equation}\label{oct1710}
\theta(x,t)=\rho(x,t)-\kappa,
\end{equation}
so that
\[
\int_\T\theta(x,t)dx= 0.
\]
Thus, the primitive functions of $\theta(x,t)$ and $G(x,t)$ are periodic.
We denote by $(\varphi, \psi)$ the mean-zero primitive functions of
$(\theta, G)$, respectively:
\begin{equation}\label{oct1724}
\theta(x,t)=\partial_x\varphi(x,t),~~\int_\T\varphi(x,t)dx=0,
\end{equation}
and
\begin{equation}\label{oct1726}
G(x,t)=\partial_x\psi(x,t),~~\int_\T\psi(x,t)dx=0.
\end{equation}
Then, $u$ can be written as
\begin{equation}\label{oct1714}
u(x,t)=\Lambda^\alpha\varphi(x,t)+\psi(x,t)+I_0(t).
\end{equation}
To determine $I_0(t)$, we use the  conservation of the momentum.
Note that the conservation law form of \eqref{eq:mainu} is
\begin{equation}\label{oct1720}
\partial_t(\rho u) +\partial_x(\rho u^2)=c_\alpha\int_\R\frac{u(t,y)-u(t,x)}{|y-x|^{1+\alpha}}
\rho(t,y)dy.
\end{equation}
Integrating (\ref{oct1720}) gives
\begin{eqnarray}\label{oct1412bis}
&&\frac{d}{dt}\int_\T\rho udx
=
c_\alpha\int_\T\int_\R\frac{u(y,t)-u(x,t)}{|y-x|^{1+\alpha}}\rho(y,t)\rho(x,t)dydx
\\
&&~~~~~~~~~~~~~~=\sum_{m\in\R\backslash\T}c_\alpha\int_\T\int_\T\frac{u(y,t)-u(x,t)}
{|y+m-x|^{1+\alpha}}\rho(y,t)\rho(x,t)dydx=0,\nonumber
\end{eqnarray}
thus
\[
%\frac{d}{dt}\int_\T \rho(x,t)u(x,t)dx=0,\quad\text{and then}\quad
\int_\T \rho(x,t)u(x,t)dx=\int_\T \rho_0(x)u_0(x)dx.
\]
Together with \eqref{oct1714}, $u$ is now uniquely defined, with $I_0(t)$ given by
\begin{equation}\label{eq:C0bis}
I_0(t)=\frac{1}{\kappa|\T|}\left[\int_\T \rho_0(x)u_0(x)dx-
\int_\T\rho(x,t)\left(\Lambda^\alpha\varphi(x,t)+\psi(x,t)\right)dx\right].
\end{equation}
Note that we have
\begin{equation}\label{oct2308}
\int_\T\rho(x,t)\Lambda^\alpha\varphi(x,t)dx=
\kappa\int_\T \Lambda^\alpha\varphi(x,t)dx+\int_\T (\partial_x\varphi(x,t))
\Lambda^\alpha\varphi(x,t)dx=0,
\end{equation}
thus
\begin{equation}\label{eq:C0}
I_0(t)=\frac{1}{\kappa|\T|}\left[\int_\T \rho_0(x)u_0(x)dx-
\int_\T\rho(x,t) \psi(x,t)dx\right].
\end{equation}
In particular, $I_0(t)$ is time-independent in the special case $G\equiv 0$, that leads to (\ref{eq:special}),
and then
we have
\begin{equation}\label{oct2310}
I_0(t)\equiv I_0(0).
\end{equation}

\subsection{The upper bound on the density}

We now prove an a priori $L^\infty$ bound on $\rho$.
\begin{thm} \label{thm:linf}
Let $\rho(x,t),u(x,t)$ be a strong solution to (\ref{eq:mainrho})-(\ref{eq:mainu})
for $0\le t\le T$,
with smooth periodic initial
conditions~$\rho_0(x)$, $u_0(x)$ such that $\rho_0(x)>0$
on $\T$.
Then, there exists a constant~$C_0>0$ that depends on~$\rho_0$ and $u_0$
but not on $T$, so that $\|\rho(\cdot,t)\|_{L^\infty}\le C_0$ for all~$t\ge 0$.
\end{thm}
This bound already indicates that  the Euler alignment system behaves not as
the fractional Burgers equation. Indeed, if we couple fractional
Burgers equation with (\ref{eq:mainrho}), the density may blow up for
$\alpha\in(0,1)$ for suitable smooth initial conditions.

%  where one may think of $(-\partial_xu)$ as the density,
% which may blow up for $\alpha\in(0,1/2)$ for suitable smooth initial conditions.

\subsubsection*{The proof of Theorem \ref{thm:linf}}

As the functions $\rho$ and $G$ obey the same continuity equation, their ratio
$F=G/\rho$ satisfies
\begin{equation}\label{oct2630}
\partial_tF+u\partial_xF=0.
\end{equation}
It follows that $F$ is uniformly bounded:
\[
\|F(\cdot,t)\|_{L^\infty}\le\|F_0\|_{L^\infty}=\left\|\frac{\partial_xu_0
-\Lambda^\alpha\rho_0}{\rho_0}\right\|_{L^\infty}<+\infty,
\]
as $\rho_0$ and $u_0$ are smooth, and $\rho_0$ is strictly positive.

In order to prove the upper bound on $\rho$,
for a fixed $t\geq0$, let $\bar{x}$ be such that
\begin{equation}\label{oct1730}
\rho(\bar{x},t)=\max_{x\in\R}\rho(x,t).
\end{equation}
It follows from \eqref{eq:rhobis} that
\begin{equation}\label{oct1802}
\partial_t\rho(\bar{x},t)=-u(\bar{x},t)\partial_x\rho(\bar{x},t)-
\rho(\bar{x},t)\partial_xu(\bar{x},t)=
-\rho(\bar{x},t)\partial_xu(\bar{x},t).
\end{equation}
Thus, to obtain an a priori upper bound on $\rho$,
it suffices to show that there exists
$C_0$ that depends on the initial conditions $\rho_0$ and $u_0$ so that
if $\rho(\bar x,t)>C_0$, then
\begin{equation}\label{oct1804}
\partial_xu(\bar{x},t)>0.
\end{equation}

To obtain (\ref{oct1804}), note that
\begin{equation}\label{oct1732}
\partial_xu=\Lambda^\alpha\rho+F\rho\geq
\Lambda^\alpha\rho-\|F_0\|_{L^\infty}\rho.
\end{equation}
%We would like to bound the right side in (\ref{oct1732}) from below,
%to show that the slope of $u$ can not become too negative.
%To this end, for a fixed $t\geq0$, let $\bar{x}$ be such that
%\begin{equation}\label{oct1730}
%\rho(\bar{x})=\max_{x\in\R}\rho(x).
%\end{equation}
In order to bound $\Lambda^\alpha\rho$ in the right side of (\ref{oct1732})
from below, we use
the nonlinear maximum principle for the fractional Laplacian,
see~\cite[Theorem 2.3]{constantin2012nonlinear}:
\begin{equation}\label{oct1736}
\text{either}\quad \Lambda^\alpha\rho(\bar{x})=\Lambda^\alpha\theta(\bar{x})\geq
\frac{\theta^{1+\alpha}(\bar{x})}{c\|\varphi\|_{L^\infty}^\alpha}\quad
\text{or}\quad \theta(\bar{x})\leq c\|\varphi\|_{L^\infty}.
\end{equation}
Here, the constant $c>0$ only depends on $\alpha$.
Recall that we denote by $\theta(x,t)$ the mean-zero
shift of~$\rho(x,t)$, as in (\ref{oct1722})
and (\ref{oct1710}), and by $\varphi(x,t)$ the mean-zero
primitive of $\theta(x,t)$, as in~(\ref{oct1724}).
Note that~$\|\varphi(\cdot,t)\|_{L^\infty}$ is uniformly
bounded:
\begin{equation}\label{oct1738}
\|\varphi(\cdot,t)\|_{L^\infty}\leq C\|\theta(\cdot,t)\|_{L^1}\leq
C\|\rho(\cdot,t)\|_{L^1}=C\|\rho_0\|_{L^1}.
\end{equation}
Therefore, if
\begin{equation}\label{oct1740}
\rho(\bar{x},t)\geq 2\kappa+C\|\rho_0\|_{L^1},
\end{equation}
with a sufficiently large $C$, which depends only on $\rho_0$ and $u_0$, then
\[
\theta(\bar x,t)=\rho(\bar x,t)-\kappa\ge 2c\|\varphi(\cdot,t)\|_{L^\infty},
\]
and the second possibility in (\ref{oct1736}) can not hold. Thus, as soon as
(\ref{oct1740}) holds, we have
\begin{equation}\label{oct1742}
\Lambda^\alpha\rho(\bar{x},t)\geq
C\frac{(\rho(\bar{x},t)-\kappa)^{1+\alpha}}{\|\rho_0\|_{L^1}^\alpha}
\geq C_0\rho(\bar{x},t)^{1+\alpha},
\end{equation}
%\[
%\Lambda^\alpha\rho(\bar{x},t)\geq
%\frac{(\rho(\bar{x},t)-\kappa)^{1+\alpha}}{c(2\kappa+c\|\rho_0\|_{L^1})^\alpha}
%\geq\frac{\rho(\bar{x},t)^{1+\alpha}}{2^{1+\alpha}
%c(2\kappa+c\|\rho_0\|_{L^1})^\alpha}
%=C\rho(\bar{x})^{1+\alpha},\]
with a constant $C_0$ that depends on the initial condition $\rho_0$.
%with a universal constant $C>0$.
Going back to (\ref{oct1732}),
this implies
\[
\partial_xu(\bar{x},t)\geq
C_0\rho(\bar{x},t)^{1+\alpha}-\|F_0\|_{L^\infty}\rho(\bar{x},t)>0.
%\quad
%\text{if }\rho(\bar{x})>\left(\frac{\|F_0\|_{L^\infty}}{C}\right)^{\frac{1}
%{\alpha}}.
\]
Thus, (\ref{oct1804}) indeed holds
if $\rho(\bar{x},t)>C_0'$, where $C_0'$ is a constant that depends only
on $\rho_0$ and $u_0$, and the proof of Theorem~\ref{thm:linf} is complete.~$\Box$
%Then, it follows from \eqref{eq:rho} that
%\[\partial_t\rho(\bar{x})=-u(\bar{x})\partial_x\rho(\bar{x})-\rho(\bar{x})\partial_xu(\bar{x})<0.\]
%Hence, $\|\rho(\cdot,t)\|_{L^\infty}=\rho(\bar{x})$ can not exceed
%\[\max\left\{2\kappa+c\|\rho_0\|_{L^1},
%  \left(\frac{\|F_0\|_{L^\infty}}{C}\right)^{\frac{1}{\alpha}}\right\},\]
%for any $t\geq0$. This provides a uniform bound on $\rho$.~$\Box$

One immediate consequence of Theorem \ref{thm:linf} is that $I_0(t)$
in \eqref{eq:C0} is uniformly bounded for all time. Indeed, it
suffices to bound
\[
\left|\int_\T\rho(x,t) \psi(x,t)dx\right|
\leq\|\rho(\cdot,t)\|_{L^2}\|\psi(\cdot,t)\|_{L^2},
\]
%\[
%\left|\int_\T\rho(x,t)\left(\Lambda^\alpha\varphi(x,t)+\psi(x,t)\right)dx\right|
%\leq\|\rho(\cdot,t)\|_{L^2}
%\left(\|\Lambda^\alpha\varphi(\cdot,t)\|_{L^2}+\|\psi(\cdot,t)\|_{L^2}\right).
%\]
%The right side is uniformly bounded due to the following
%estimates: first, we simply have
%\begin{equation}\label{oct1810}
%\|\rho(\cdot,t)\|_{L^2}\leq C\|\rho(\cdot,t)\|_{L^\infty}.
%\end{equation}
%As $\alpha\in(0,1)$, relation (\ref{oct1724}) implies that
%the term $\Lambda^\alpha\varphi$ can be estimated as
%\begin{equation}\label{oct1812}
%\|\Lambda^\alpha\varphi(\cdot,t)\|_{L^2}\leq\|\theta(\cdot,t)\|_{L^2}
%\leq C\|\rho(\cdot,t)\|_{L^\infty}.
%\end{equation}
%Finally, we note that
while
\begin{equation}\label{oct1814}
\|\psi(\cdot,t)\|_{L^2}\leq C\|G(\cdot,t)\|_{L^2}
\leq C\|G(\cdot,t)\|_{L^\infty}
\leq C\|\rho(\cdot,t)\|_{L^\infty}\|F_0\|_{L^\infty}\leq C,
\end{equation}
where $C$ is a universal constant independent of $t$. Summarizing, we have
\begin{equation}\label{oct1806}
|I_0(t)|\le C_0,
\end{equation}
with a constant $C_0$ that depends only on $\rho_0$ and $u_0$.

Thus, we have the following a priori bound on $\|u\|_{L^2}$.
\begin{cor} \label{cor:ul2}
Let $\rho(x,t),u(x,t)$ be
a strong solution to (\ref{eq:mainrho})-(\ref{eq:mainu})
for $0\le t\le T$,
with smooth periodic initial
conditions~$\rho_0(x)$, $u_0(x)$ such that $\rho_0(x)>0$
on $\T$.
There exists a constant~$C_0$ that depends only on $\rho_0$ and $u_0$
but not not on $T$ so that
$\|u(\cdot,t)\|_{L^2}\le C_0$ for all~$0\le t\le T$.
\end{cor}
{\bf Proof.} This follows immediately from the bound
\[
\|u(\cdot,t)\|_{L^2}\leq \|\Lambda^\alpha\varphi(\cdot,t)\|_{L^2}
+\|\psi(\cdot,t)\|_{L^2}+|I_0(t)|,
\]
together with the bound
\begin{equation}\label{oct1812}
\|\Lambda^\alpha\varphi(\cdot,t)\|_{L^2}\leq C\|\theta(\cdot,t)\|_{L^2}
\leq C\|\rho(\cdot,t)\|_{L^\infty},
\end{equation}
and (\ref{oct1814})-(\ref{oct1806}).~$\Box$

The uniform upper bound on the density also implies a uniformly Lipschitz bound on $F$.
\begin{lem}\label{lem:oct2604}
The function $F=G/\rho$ is Lipschitz, and the Lipschitz bound is uniform in time.
\end{lem}
{\bf Proof.}
Recall that $F$ satisfies (\ref{oct2630}), thus $p=\partial_xF$ satisfies the same continuity equation as $\rho$:
\begin{equation}\label{oct2632bis}
\partial_tp+\partial_x(up)=0,
\end{equation}
%So $\partial_xF$ satisfies the same continuity equation as $\rho$.
and $w=p/\rho$ is a solution of
\[
\partial_tw+u\partial_xw=0.
\]
It follows that $\|w(\cdot,t)\|_{L^\infty}=\|w_0\|_{L^\infty}$, and
therefore,
\[
\|\partial_xF(\cdot,t)\|_{L^\infty}\leq \|w_0\|_{L^\infty}\|\rho(\cdot,t)\|_{L^\infty}.
\]
Theorem \ref{thm:linf} implies now that $F$ is Lipschitz, with
a time-independent Lipschitz bound.~$\Box$

\subsection{A lower bound on the density}
A uniform lower bound on $\rho$ plays an important role as it keeps
the dissipation active. The following lemma ensures no creation of
vacuum in finite time.

\begin{lem}%[Lower bound on $\rho$]
\label{lem:lower}
Let $\rho(x,t),u(x,t)$ be
a strong solution to (\ref{eq:mainrho})-(\ref{eq:mainu})
for $0\le t\le T$,
with smooth periodic initial
conditions~$\rho_0(x)$, $u_0(x)$ such that $\rho_0(x)>0$
on $\T$. There
exists a positive constant~$C_0>0$ that depends on $\rho_0$ and $u_0$ but
not on $T$, so that
\begin{equation}\label{oct1820}
\rho(x,t)\geq\farc{1}{C_0(1+t)},
\quad \hbox{ for all $x\in\T$ and $0\le t\le T$.}
\end{equation}
%decreasing function $\rho_m(t)$ that does
%not depend on $T$, such that
%\begin{equation}\label{oct1820}
%\rho(x,t)\geq\rho_m(t),\quad \hbox{ for all $x\in\T$ and $0\le t\le T$.}
%\end{equation}
\end{lem}
{\bf Proof.}
Fix some $t>0$ and let $\underline{x}$ be such that
\[
\rho(\underline{x},t)=\min_x\rho(x,t).
\]
Then we have
\[
\Lambda^\alpha\rho(\underline{x},t) \leq 0,
\]
and thus
\begin{equation}\label{oct2814}
\rho_m(t)=\rho(\underline x,t)=\min_{x\in\T}\rho(x,t),
\end{equation}
satisfies
\begin{eqnarray*}
&&\farc{d\rho_m(t)}{dt}=\partial_t\rho(\underline{x},t)=
[-\partial_xu(\underline{x},t)]\rho(\underline{x},t)
\geq-\big(\Lambda^\alpha\rho(\underline{x},t)
+\|F_0\|_{L^\infty}\rho_m(t)\big)\rho_m(t)\\
&&~~~~~~~~~\geq -\|F_0\|_{L^\infty}\rho_m(t)^2.
\end{eqnarray*}
%where we use the fact $\Lambda^\alpha\rho(\underline{x})<0$ for the
%last inequality.
If the minimum is achieved at more than one point, we just need to take a minimum over all of them in the above estimate,
which leads to the same bound.
Notice that $\rho_m(t)$ is Lipschitz in time, so the estimate is valid for a.e. $t$, and ${d\rho_m}/{dt}$
determines $\rho_m(t).$
Integrating this differential inequality, we get
\begin{equation}\label{oct2810}
\rho_m(t)\geq\frac{1}{[\rho_m(0)]^{-1}+t\|F_0\|_{L^\infty}},
\end{equation}
finishing the proof.~$\Box$
%From the proof, we know that the lower bound $\rho_m$ decays to zero
%like $t^{-1}$, when $t\to\infty$.

In particular, in the special case $G\equiv 0$, that is, for (\ref{eq:special}) we have the following.
\begin{cor}\label{cor:oct2502}
Let $\rho(x,t)$ be the solution of (\ref{eq:special}). Then, we have
\begin{equation}\label{oct2502}
\rho(x,t)\ge \min_{x\in\T}\rho_0(x),\hbox{ for all $t>0$ and $x\in\T$.}
\end{equation}
\end{cor}

\section{The local wellposedness}\label{sec:3}
% Let us consider the case that  $x\in\mathbb{T}$.
% Given smooth initial data $u_0$, we pick $\rho_0$ such that
% $\rho_0(x)\geq\rho_m>0$, and
% $\Lambda^\alpha\rho_0(x)=\partial_xu_0(x)$. Under this setup,
% $C_3\equiv0$, and therefore $\partial_xu(x,t)=\Lambda^\alpha\rho(x,t)$ for
% all $(x,t)\in\mathbb{T}\times\mathbb{R}^+$. The dynamics becomes
% \begin{equation}\label{eq:special}
% \partial_t\rho+\partial_x(\rho
% u)=0,\quad \partial_xu=\Lambda^\alpha\rho.
% \end{equation}
% In particular, for $\alpha=1$, the equation is known as a 1D model of
% the quasi-geostrophic equation, and has been studied in \cite{chae2005finite},
% where global regularity is obtained. It indicates that there is no
% finite time blowup for \eqref{eq:main}.
% Our main interest is the case $\alpha<1$, where fractal Burgers
% equation \eqref{eq:fractalBurgers} blows up in finite time.

The a priori bounds on $\rho$ established in the previous section rule out some kinds of finite time blow up,
but do not imply that there is no finite time
shock formation. This remains to be shown. To proceed further, we first establish a local well-posedness theory
for solutions of the Euler alignment system with smooth initial conditions.

\begin{thm}\label{thm:local}
Let $\alpha \in (0,1).$ Assume that the initial conditions $\rho_0$ and $u_0$ satisfy
\begin{equation}\label{oct1830}
\rho_0\in H^s(\T),\quad \min_{x\in\T}\rho_0(x)>0,\quad
\partial_xu_0-\Lambda^\alpha\rho_0\in H^{s-\frac{\alpha}{2}}(\T),
\end{equation}
with a sufficiently large even integer $s>0$. Then, there exists $T_0>0$
such that the system~(\ref{eq:mainrho})-(\ref{eq:mainu})
has a unique strong solution $\rho(x,t), u(x,t)$ on $[0,T_0]$, with
\begin{equation}\label{oct1832}
\rho\in C([0,T_0], H^s(\T))\times L^2([0,T_0], H^{s+\frac{\alpha}{2}}(\T)),\quad
u\in C([0,T_0], H^{s+1-\alpha}(\T)).
\end{equation}
Moreover, a necessary and sufficient condition for the solution to
exist on a time interval~$[0,T]$ is
\begin{equation}\label{eq:BKM}
\int_0^T\|\partial_x\rho(\cdot,t)\|_{L^\infty}^2 dt<\infty.
\end{equation}
\end{thm}

Condition \eqref{eq:BKM} is a Beale-Kato-Majda type
criterion. It indicates that the solution is globally regular if
$\partial_x\rho$ is uniformly bounded in the $L^\infty$ norm. We will show
that such bound actually does hold in Section~\ref{sec:4},
using the modulus of continuity method.

\subsection{The commutator estimates}

We will need some commutator estimates for the local well-posedness theory.
We will use the following notation:
\begin{eqnarray*}
&&[\mathcal{L}, f, g]=\mathcal{L} (fg)-f \mathcal{L} g-g \mathcal{L} f,\\
&&[\mathcal{L}, f]g=\mathcal{L} (fg)-f \mathcal{L} g.
\end{eqnarray*}
\begin{lem}
The following commutator estimates hold:\\
(i) for any $n\ge 1$,  we have
\begin{equation}\label{est:s}
\|[\partial_x^n, f, g]\|_{L^2}\le C \big(
\|\partial_xf\|_{L^\infty}\|g\|_{H^{n-1}}+
\|\partial_xg\|_{L^\infty}\|f\|_{H^{n-1}}\big),%\quad n\ge 1,
\end{equation}
(ii) for any $\gamma\in(0,1)$ and $\epsilon>0$, we have
\begin{equation}\label{est:alpha}
\|[\Lambda^\gamma, f, g]\|_{L^2}\le C
\|f\|_{L^2}\|g\|_{C^{\gamma+\epsilon}} ,
%\quad \gamma\in(0,1),~\epsilon>0
\end{equation}
(iii) for any $\gamma>0$, we have
\begin{equation}\label{est:2}
\|[\Lambda^\gamma, f]g\|_{L^2}\le C \big(
\|\partial_xf\|_{L^\infty}\|g\|_{H^{\gamma-1}}+\|f\|_{H^\gamma}\|g\|_{L^\infty}
\big).
\end{equation}
\end{lem}
Let us comment briefly on the proof of these estimates.
Estimate \eqref{est:s} can be obtained by
the standard Gagliardo-Nirenberg interpolation inequality. As
$\Lambda^2=-\partial^2_{xx}$, this estimate holds if we replace the
operator $\partial_x^n$ by $\Lambda^s$ with an even integer $s$.

A version of \eqref{est:alpha} is discussed in \cite[Theorem
A.8]{kenig1993well}.
We sketch the proof in Appendix~\ref{sec:app}.
Finally, estimate \eqref{est:2} is due to Kato and Ponce
\cite{kato1988commutator}. The proof is similar to that of (\ref{est:alpha}).
%and $\Lambda^\gamma$
%can be replaced by $H\Lambda^\gamma$, where $H$ is the Hilbert transform.

\subsection{The proof of the local well-posedness}

%We are ready to establish a local wellposedness theory for
%(\ref{eq:mainrho})-(\ref{eq:mainu}).
It will be convenient to use the variables
$(\theta, G)$, so that equations (\ref{eq:rhobis})-(\ref{eq:Gbis})
take the form
\begin{eqnarray}\label{oct1828}
&&\partial_t\theta+\partial_x(\theta u)=-\kappa \partial_xu,\quad
\partial_t G+\partial_x (G u)=0,\\
&&\partial_xu=\Lambda^\alpha\theta+G.\label{oct2002}
\end{eqnarray}
Here $\kappa$ is the constant in time mean of $\rho$, as in (\ref{oct1722}).

Let us fix $T>0$ and take a sufficiently large even integer $s>0$.
We will aim to obtain a differential inequality
on
\begin{equation}\label{oct2004}
Y(t):=1+\|\theta(\cdot,t)\|_{H^s}^2+\|G(\cdot,t)\|_{H^{s-\frac{\alpha}{2}}}^2,
\end{equation}
that will have bounded solutions on a time interval $[0,T_0]$, with a sufficiently small $T_0$ depending
on the initial conditions.
To this end, we apply the operator
$\Lambda^s$ to the   equation for~$\theta$ in (\ref{oct1828}), multiply the result
by~$\Lambda^s\theta$ and integrate in~$x$:
\begin{equation}\label{eq:Hstheta}
\frac{1}{2}\frac{d}{dt}\|\theta(\cdot,t)\|_{\dot{H}^s}^2=
-\int\big(\Lambda^s\theta\cdot\Lambda^s\partial_x(\theta u)\big)dx
-\kappa\|\theta(\cdot,t)\|_{\dot{H}^{s+\frac{\alpha}{2}}}^2
-\kappa\int\big(\Lambda^s\theta\cdot\Lambda^sG\big) dx.
\end{equation}
The second term in the right side produces
the dissipation. We shall use it to control the
other two terms.

We split the first term in the right side
of \eqref{eq:Hstheta} into three pieces:
\begin{align}
\int\Lambda^s\theta\cdot\Lambda^s\partial_x(\theta u)dx
=&\int(\Lambda^s\theta\cdot\Lambda^s\partial_xu)\theta dx+
\int(\Lambda^s\theta\cdot u)(\Lambda^s\partial_x\theta) dx+
\int\Lambda^s\theta\cdot[\Lambda^s\partial_x , u, \theta]dx\nonumber\\
=&~ I+II+III.\label{oct1920}
\end{align}
Let us start with $I$: % It can be decomposed into three parts.
\begin{align}\label{oct1910}
I=&\int(\Lambda^{s-\frac{\alpha}{2}}\partial_xu)\cdot
\Lambda^{\frac{\alpha}{2}}(\theta\cdot\Lambda^s\theta)dx\nonumber\\
=&\int(\Lambda^{s-\frac{\alpha}{2}}\partial_xu)
\cdot(\Lambda^{s+\frac{\alpha}{2}}\theta)\cdot\theta dx+
\int(\Lambda^{s-\frac{\alpha}{2}}\partial_xu)\cdot (\Lambda^s\theta)\cdot
(\Lambda^{\frac{\alpha}{2}}\theta) dx+
\int(\Lambda^{s-\frac{\alpha}{2}}\partial_xu)\cdot
[\Lambda^{\frac{\alpha}{2}}, \Lambda^s\theta, \theta]dx\nonumber\\
=&~I_1 +I_2+I_3.
\end{align}
For $I_1$, we have, using (\ref{oct2002}):
\[
I_1=\int|\Lambda^{s+\frac{\alpha}{2}}\theta|^2\cdot\theta dx
+\int(\Lambda^{s-\frac{\alpha}{2}}
G)\cdot(\Lambda^{s+\frac{\alpha}{2}}\theta)\cdot\theta
dx= I_{11} +I_{12}.
\]
The term $I_{11}$ is controlled by the dissipation in the right side
of (\ref{eq:Hstheta}): set
\[
\rho_m(t)=\inf_{0\le \tau\le t,x\in\T}\rho(x,\tau).
\]
Note that $\rho_m(t)>0$ by Lemma~\ref{lem:lower}. Then we have, using
Lemma \ref{lem:lower}:
\begin{equation}\label{oct1912}
-I_{11}-\kappa\|\theta\|_{\dot{H}^{s+\frac{\alpha}{2}}}^2\leq
(\|\theta_-\|_{L^\infty}-\kappa)\|\theta\|_{\dot{H}^{s+\frac{\alpha}{2}}}^2
\leq -\rho_m(t)\|\theta\|_{\dot{H}^{s+\frac{\alpha}{2}}}^2.
\end{equation}
To bound $I_{12}$ we use the
H\"{o}lder inequality:
\begin{equation}\label{oct1914}
|I_{12}|\leq\|G\|_{\dot{H}^{s-\frac{\alpha}{2}}}\|\theta\|_{\dot{H}^{s+\frac{\alpha}{2}}}
\|\theta\|_{L^\infty}
\leq\frac{\rho_m}{6}\|\theta\|_{\dot{H}^{s+\frac{\alpha}{2}}}^2+
\frac{3}{2\rho_m}\|\theta\|_{L^\infty}^2\|G\|_{\dot{H}^{s-\frac{\alpha}{2}}}^2.
\end{equation}
In order to control the term $I_2$ in (\ref{oct1910}), we, once again,
use (\ref{oct2002}), and
the H\"{o}lder inequality:
\begin{align}
|I_2|\leq&\left(\|\theta\|_{\dot{H}^{s+\frac{\alpha}{2}}}+
\|G\|_{\dot{H}^{s-\frac{\alpha}{2}}}\right)\|\theta\|_{\dot{H}^s}
\|\Lambda^{\frac{\alpha}{2}}\theta\|_{L^\infty}\nonumber\\
\leq&~\frac{\rho_m}{6}\|\theta\|_{\dot{H}^{s+\frac{\alpha}{2}}}^2+
\left(\frac{3}{2\rho_m}+\frac{1}{2}\right)\|\Lambda^{\frac{\alpha}{2}}\theta\|_{L^\infty}^2\|\theta\|_{\dot{H}^s}^2
+\frac{1}{2}\|G\|_{\dot{H}^{s-\frac{\alpha}{2}}}^2.\label{oct1916}
\end{align}
The contribution of
$I_3$ in (\ref{oct1910})
is bounded using the commutator estimate \eqref{est:alpha}:
\begin{align}
|I_3|\leq&\left(\|\theta\|_{\dot{H}^{s+\frac{\alpha}{2}}}+
\|G\|_{\dot{H}^{s-\frac{\alpha}{2}}}\right)
\|[\Lambda^{\frac{\alpha}{2}}, \Lambda^s\theta, \theta]\|_{L^2}
\leq C\left(\|\theta\|_{\dot{H}^{s+\frac{\alpha}{2}}}+
\|G\|_{\dot{H}^{s-\frac{\alpha}{2}}}\right)\|\theta\|_{\dot{H}^s}
\|\theta\|_{C^{\frac{\alpha}{2}+\epsilon}}\nonumber\\
\leq&~
\frac{\rho_m}{6}\|\theta\|_{\dot{H}^{s+\frac{\alpha}{2}}}^2+
\left(\frac{3}{2\rho_m}+\frac{1}{2}\right)C^2\|\theta\|_{C^{\frac{\alpha}{2}+\epsilon}}^2\|\theta\|_{\dot{H}^s}^2
+\frac{1}{2}\|G\|_{\dot{H}^{s-\frac{\alpha}{2}}}^2.\label{oct1922}
\end{align}
Next, we estimate the term $II$ in (\ref{oct1920}),
integrating by parts
\begin{equation}\label{oct1924}
|II|=\frac{1}{2}\left|\int(\Lambda^s\theta)^2
\cdot \partial_xu~dx\right|\le C
\left(\|\Lambda^\alpha\theta\|_{L^\infty}+\|G\|_{L^\infty}\right)\|\theta\|_{\dot{H}^s}^2.
\end{equation}
For the term $III$ in (\ref{oct1920}),
we apply the commutator estimate \eqref{est:s} and
get
\begin{equation}\label{oct1934}
|III|\leq\|\theta\|_{\dot{H}^s}\|[\Lambda^s\partial_x , u, \theta]\|_{L^2}
\le C\|\theta\|_{H^s}\left(\|\partial_xu\|_{L^\infty}\|\theta\|_{H^s}
+\|\partial_x\theta\|_{L^\infty}\|u\|_{H^s}\right).
%\nonumber\\
%\le &C(\|\partial_x\theta\|_{L^\infty}+\|G\|_{L^\infty})
%(\|\theta\|_{H^s}^2+\|G\|_{H^{s-1}}^2+1).\label{oct1930}
\end{equation}
To estimate $\|u\|_{H^s}$ in the right side,
we apply Corollary \ref{cor:ul2} to get
\begin{equation}\label{oct1936}
\|u\|_{H^s}=\|u\|_{L^2}+\|\partial_xu\|_{H^{s-1}}
\leq C(1+\|\theta\|_{H^{s-1+\alpha}}+\|G\|_{H^{s-1}}).
\end{equation}
We also have, using the uniform bound on the density:
\begin{equation}\label{oct2006}
|\Lambda^\alpha\theta|\le c_\alpha\int_\R
\farc{|\theta(x)-\theta(y)|dy}{|x-y|^{1+\alpha}}\le
C(\|\theta\|_{L^\infty}+\|\partial_x\theta\|_{L^\infty})\le
C_0(1+\|\partial_x\theta\|_{L^\infty}),
\end{equation}
with a constant $C_0$ that depends on $\rho_0$ and $u_0$.
Therefore, $\partial_xu$ satisfies
\begin{equation}\label{oct1938}
\|\partial_xu\|_{L^\infty}\le \|\Lambda^\alpha\theta\|_{L^\infty}+\|G\|_{L^\infty}
\le C(1+\|\partial_x\theta\|_{L^\infty}+\|G\|_{L^\infty}).
\end{equation}
Together, (\ref{oct1934})-(\ref{oct1938}) give
%\begin{equation}\label{oct2006}
%|III|\le C(1+\|\partial_x\theta\|_{L^\infty}+\|G\|_{L^\infty})
%\big[\|\theta\|_{H^s}^2+
%(1+\|\theta\|_{H^{s-1+\alpha}}+\|G\|_{H^{s-1}})\|\theta\|_{H^s}\big]
%%(1+\|\theta\|_{H^s}^2+\|G\|_{H^{s-1}}^2)
%\end{equation}
%
%\begin{equation}\label{oct1932}
%|III|\le C\|\theta\|_{H^s}(\|\partial_x\theta\|_{L^\infty}+\|G\|_{L^\infty})
%(1+\|\theta\|_{H^s}^2+\|G\|_{H^{s-1}}^2)
%\end{equation}
\begin{equation}\label{oct1932bis}
|III|\le C(1+\|\partial_x\theta\|_{L^\infty}+\|G\|_{L^\infty})
(1+\|\theta\|_{H^s}^2+\|G\|_{H^{s-1}}^2).
\end{equation}
The third term in the right side of~\eqref{eq:Hstheta} can be estimated as
\[
\kappa\left|\int(\Lambda^s\theta)\cdot(\Lambda^sG )dx\right|\leq
\kappa\|\theta\|_{\dot{H}^{s+\frac{\alpha}{2}}}\|G\|_{\dot{H}^{s-\frac{\alpha}{2}}}
\leq \frac{\rho_m}{6} \|\theta\|_{\dot{H}^{s+\frac{\alpha}{2}}}^2+
\frac{3\kappa^2}{2\rho_m}\|G\|_{\dot{H}^{s-\frac{\alpha}{2}}}^2.
\]
%From Theorem \ref{thm:linf}, we know $\|\theta\|_{L^\infty}$ and
%$\|G\|_{L^\infty}$ are uniformly bounded. So we absorb them into the
%universal constant $C$. Also, it is easy to show
%\[\|\Lambda^\gamma\theta\|_{L^\infty}\lesssim\|\partial_x\theta\|_{L^\infty},
%\quad\forall~\gamma\in(0,1).\]
Putting the above estimates together,
we end up with the following inequality:
\begin{equation}\label{oct2010}
\frac{1}{2}\frac{d}{dt}\|\theta\|_{H^s}^2\leq
C\left(1+\frac{1}{\rho_m}\right)
(1+\|\partial_x\theta\|_{L^\infty}^2+\|G\|_{L^\infty})
(\|\theta\|_{H^s}^2+\|G\|_{H^{s-\frac{\alpha}{2}}}^2+1)-
\frac{\rho_m}{3}\|\theta\|_{H^{s+\frac{\alpha}{2}}}^2.
\end{equation}
In order to close the estimate, and obtain a bound on $Y(t)$ defined in (\ref{oct2004}),
we write:
\begin{align}\label{oct2020}
\frac{1}{2}\frac{d}{dt}\|G\|_{\dot{H}^{s-\frac{\alpha}{2}}}^2
=&~-\int(\Lambda^{s-\frac{\alpha}{2}} G)\cdot(\Lambda^{s-\frac{\alpha}{2}}\partial_x (Gu))dx\\
=&~-\int(\Lambda^{s-\frac{\alpha}{2}}
   G)\cdot(u \Lambda^{s-\frac{\alpha}{2}}\partial_x G)~dx-
\int(\Lambda^{s-\frac{\alpha}{2}}G)\cdot[\Lambda^{s-\frac{\alpha}{2}}\partial_x, u]G dx=
IV+V.\nonumber
\end{align}
The term $IV$ can be treated as $II$ via integration by parts, together with (\ref{oct2006}):
\begin{equation}\label{oct2016}
|IV|=\frac{1}{2}\left|\int(\Lambda^{s-\frac{\alpha}{2}} G)^2
\cdot \partial_xu~dx\right|\leq
C\left(1+\|\partial_x\theta\|_{L^\infty}+\|G\|_{L^\infty}\right)
\|G\|_{\dot{H}^{s-\frac{\alpha}{2}}}^2.\
\end{equation}
To bound $V$, we apply the commutator estimate \eqref{est:2}, as well as (\ref{oct1938}):
\begin{align}\label{oct2014}
|V|\leq&~\|G\|_{\dot{H}^{s-\frac{\alpha}{2}}}
\|[\Lambda^{s-\frac{\alpha}{2}}\partial_x, u]G\|_{L^2}
\leq C\|G\|_{\dot{H}^{s-\frac{\alpha}{2}}}\left(\|\partial_xu\|_{L^\infty}
\|G\|_{H^{s-\frac{\alpha}{2}}}
+\|G\|_{L^\infty}\|\partial_xu\|_{H^{s-\frac{\alpha}{2}}}\right)\nonumber\\
\leq&~C (1+\|\partial_x\theta\|_{L^\infty}+\|G\|_{L^\infty}) \|G\|_{H^{s-\frac{\alpha}{2}}}^2+
C\|G\|_{L^\infty}\|G\|_{H^{s-\frac{\alpha}{2}}}\|\theta\|_{H^{s+\frac{\alpha}{2}}}\nonumber\\
\leq&~ \frac{\rho_m}{6}\|\theta\|_{H^{s+\frac{\alpha}{2}}}^2 +C
\left(1+\frac{1}{\rho_m}\|G\|_{L^\infty}^2+\|\partial_x\theta\|_{L^\infty}+\|G\|_{L^\infty})\right) \|G\|_{H^{s-\frac{\alpha}{2}}}^2.
\end{align}
%
%Define
%$Y_s(t):=\|\theta(\cdot,t)\|_{H^s}^2+\|G(\cdot,t)\|_{H^{s-\frac{\alpha}{2}}}^2+C$. Our
Now, estimates (\ref{oct2010})-(\ref{oct2014}), together with the uniform bound on $\|G\|_{L^\infty}$, yield an inequality
\begin{equation}\label{eq:Ys}
\frac{d}{dt}Y(t)\leq C\left(1+\frac{1}{\rho_m(t)}\right)
(1+\|\partial_x\theta(\cdot,t)\|_{L^\infty}^2)Y(t)
-\frac{\rho_m(t)}{6}\|\theta(\cdot,t)\|_{H^{s+\frac{\alpha}{2}}}^2.
\end{equation}
For $s>{3}/2$, $H^s$ is embedded in $W^{1,\infty}$. This,
together with Lemma \ref{lem:lower}, implies
\begin{equation}\label{oct2032}
\frac{d}{dt}Y(t)\leq C(1+t)(1+Y(t))Y(t),
\end{equation}
and the local in time well-posedness for solutions  with $H^s$ initial data follows. Moreover, it follows
from (\ref{eq:Ys}) that
\begin{equation}\label{oct2034}
Y(T)\leq
Y(0)\exp\left[C\int_0^T(1+t)(1+\|\partial_x\theta(\cdot,t)\|_{L^\infty}^2)dt\right].
\end{equation}
For all finite $T>0$, if the Beale-Kato-Majda criterion \eqref{eq:BKM}
is satisfied, the right   side of (\ref{oct2034})
is finite, whence
\[
\theta\in C([0,T],H^s(\T)),\quad
G(\cdot,t)\in C([0,T],H^{s-\frac{\alpha}{2}}(\T)),
\]
and thus $\rho\in C([0,T],H^s(\T))$. Furthermore, integrating \eqref{eq:Ys} in
$[0,T]$, we see that if (\ref{eq:BKM}) holds, then
\begin{align*}
\frac{\rho_m(T)}{6}&\|\theta\|_{L^2([0,T],H^{s+\frac{\alpha}{2}})}^2<+\infty,
%\leq
%\int_0^T\frac{\rho_m(t)}{6}\|\theta(\cdot,t)\|_{H^{s+\frac{\alpha}{2}}}^2dt\\
%&\leq Y_s(0)-Y_s(T)+C\left(1+\frac{1}{\rho_m(T)}\right)
%\int_0^T(1+\|\partial_x\theta(\cdot,t)\|_{L^\infty}^2)Y_s(t)dt<\infty.
\end{align*}
%It implies $\theta$, as well as
thus $\rho\in L^2([0,T],H^{s+\frac{\alpha}{2}})$.
To recover the conditions on $u$ in (\ref{oct1832}), we apply Corollary
\ref{cor:ul2} and get
\[\|u(\cdot,t)\|^2_{H^{s+1-\alpha}}=\|u(\cdot,t)\|^2_{L^2}
+\|\partial_xu(\cdot,t)\|^2_{H^{s-\alpha}}
\leq C+CY(t)<\infty.
\]
%The continuity in time is easy to check as $C_0(t)$ is continuous if
%$Y_s(t)$ is bounded.
This ends the proof
of Theorem \ref{thm:local}.

\section{The global regularity}\label{sec:4}

In this section, we derive a uniform~$L^\infty$-bound on
$\partial_x\rho$, using a variant of the modulus of continuity method.
Together with the Beale-Kato-Majda type criterion \eqref{eq:BKM},
this will imply the global well-posedness of the Euler alignment
system \eqref{eq:mainrho}-\eqref{eq:mainu}, and prove
Theorem~\ref{thm:global}. We will first consider the special case $G\equiv 0$,
that is, the system (\ref{eq:special}). The nonlinear diffusive term makes the 
problem subtler than in the SQG or Burgers equation case. 
Finally, we prove the result to the general Euler alignment system, 
using a combination of an appropriate scaling argument, estimate on the minimum of $\rho,$ 
and additional regularity estimates. In this case, the bound on $\partial_x \rho$ 
will depend on time and may grow, but remains finite for every $t >0.$ 

For convenience, we work on~$\R$, and extend $\rho$ and $u$ periodically in
space.

\subsection{The modulus of continuity}

We say that a function $f$ obeys modulus of continuity $\omega$ if
\[
f(x)-f(y)<\omega(|x-y|),\quad\text{for all }x,y\in\R.
\]
We will work with the following modulus of continuity for the density $\rho$:
\begin{equation}\label{oct2102}
\omega(\xi)=\begin{cases}\xi-\xi^{1+{\alpha}/{2}},&0\le\xi<\delta\\
\gamma\log(\xi/\delta)+\delta-\delta^{1+\alpha/2},&\xi\geq\delta,
\end{cases}
\end{equation}
%as well as its rescalings
%\begin{equation}\label{oct2302}
%\omega_B(\xi)=\omega(B\xi),
%\end{equation}
%with $B>0$.
%The constant $c$ in (\ref{oct2102}) is chosen as
%\begin{equation}\label{oct2416}
%c=\delta-\delta^{1+\frac{\alpha}{2}}-\gamma\log\delta,
%\end{equation}
so that $\omega$
is continuous at $\xi=\delta$. The parameters $\delta$ and $\gamma$
are sufficiently small positive numbers to be specified later.
%We then generate a family of moduli of continuity $\omega_B$ by
%\[
%\omega_B(\xi)=\omega(B\xi).
%\]
The modulus $\omega$ is continuous, piecewise differentiable, increasing and concave,  and satisfies
\begin{equation}\label{oct2104}
\hbox{ %$\omega_B(0+)=0$,
$\omega''(0)=-\infty$.}
\end{equation}
%and $\omega_B'(0+)=B$.
The following proposition describes the only
possible modulus breakthrough scenario for evolution equations.
\begin{prop}[\cite{kiselev2007global}]
Suppose $\rho_0$ obeys a modulus of continuity $\omega$ that satisfies (\ref{oct2104}).
If the solution~$\rho(x,t)$ violates
$\omega$ at some positive time, then there must exist $t_1>0$ and
$x_1\neq y_1$ such that
\begin{equation}\label{oct2110}
\rho(x_1,t_1)-\rho(y_1,t_1)=\omega(|x_1-y_1|),
\hbox{  and $\rho(\cdot,t)$ obeys $\omega$ for every $0\leq t<t_1$.}
\end{equation}
\end{prop}
Thus, to prove that $\rho$ obeys a modulus of continuity $\omega$ for all times $t>0$, it is sufficient to
prove that if (\ref{oct2110}) holds, then
\begin{equation}\label{oct2112}
\partial_t(\rho(x_1,t_1)-\rho(y_1,t_1))<0.
\end{equation}
As a remark on the notation,
we will again use $C$ as a notation for various
universal constants that do not
depend on $T, \delta$ and $\gamma$.

\subsection{The global regularity for the special system with $G\equiv0$}

Let us first consider the special case $G\equiv 0$,
or, equivalently, the system \eqref{eq:special}:
%This corresponds to the special system \eqref{eq:special}:
\begin{equation}\label{oct2204}
\partial_t\rho+\partial_x(\rho
u)=0,\quad \partial_xu=\Lambda^\alpha\rho.
\end{equation}
As the mean of $u$ is preserved by the evolution -- see (\ref{oct2310}),
we may assume without loss of generality that
\begin{equation}\label{oct2312}
\int_\T u(x,t)dx=0,
\end{equation}
for otherwise we would simply consider (\ref{oct2204}) in a frame moving
the speed equal to the mean of $u_0$. Thus, we have
\begin{equation}\label{oct2342}
u(x,t)=\Lambda^\alpha\varphi(x,t).
\end{equation}
Here, $\varphi(x,t)$ is the mean-zero primitive of $\theta(x,t)$, as in
(\ref{oct1724}).
%In other words, the term $u_2$ in (\ref{oct2202}) vanishes identically,
%and we only need to analyze $u_1$.
We will prove the following result.
\begin{thm}\label{thm:special}
The system \eqref{oct2204} with a smooth periodic
initial condition $\rho_0$ such that~$\rho_0(x)>0$ for all $x\in\T$ has a unique global
smooth solution.
\end{thm}
The key step in the proof is
\begin{lem}\label{lem:oct2402}
Suppose that $m={\rm min}_{x \in \T} \rho_0(x) >0.$ Then there exist $\delta_m$ and $\gamma_m,$ independent of the period
of the initial data, such that if $\rho_0(x)$ obeys the modulus of continuity
$\omega$ given by~(\ref{oct2102}), then  $\rho(x,t)$ obeys $\omega$
for all~$t>0$.
\end{lem}
Theorem~\ref{thm:special} is a consequence of Lemma~\ref{lem:oct2402}. Indeed, suppose that Lemma~\ref{lem:oct2402} is true.
Notice that the equation \eqref{oct2204} has a scaling invariance:
if $\rho(x,t)$ is a solution, then so is
\begin{equation}\label{oct2504}
\rho_\lambda(x,t)=\rho(\lambda
x,\lambda^\alpha t),
\end{equation}
for any $\lambda>0$.
% For any initial condition $\rho_0(x)>0$ let
% \begin{equation}\label{oct2506}
% m=\min_{x\in\T}\rho_0(x).
% \end{equation}
From the properties of the modulus of continuity $\omega$ given by (\ref{oct2102}) (in particular its growth at infinity) it follows that we
can find $\lambda>0$ sufficiently small such  that $\rho_\lambda^0(x)=\rho_0(\lambda x)$ obeys
$\omega$ with $\delta=\delta_m$, $\gamma=\gamma_m$ provided by Lemma~\ref{lem:oct2402}.
Note that the rescaling (\ref{oct2504}) does not change the minimum
of $\rho$. As $\delta_m$ and $\gamma_m$ do not depend on the period,
Lemma~\ref{lem:oct2402} shows that $\rho_\lambda(x,t)$
obeys $\omega$ for all $t>0$. In particular, it follows that
\begin{equation}\label{oct2606}
|\partial_x\rho_\lambda(t,x)|\le 1,~~\hbox{ for all $t>0$ and $x\in\T$.}
\end{equation}
As we have mentioned, (\ref{oct2606}) together with the Beale-Kato-Majda type criterion \eqref{eq:BKM},
implies that~$\rho_\lambda(t,x)$
is a global in time
solution of (\ref{oct2204}), and thus so is $\rho(t,x)$.
%which is
%Finally, as $\rho$ obeys $\omega$ for all $t>0$, we know that
%\begin{equation}\label{oct2606}
%|\partial_x\rho(t,x)|\le 1,~~\hbox{ for all $t>0$ and $x\in\T$.}
%\end{equation}
%As we have mentioned, (\ref{oct2606}) together with the Beale-Kato-Majda type criterion \eqref{eq:BKM},
%implies the conclusion of Theorem~\ref{thm:special}.

%\begin{proof}
% For any fixed $T>0$, there exists a family of modulus $\omega_B$
% which can not be broken through.
%
%Given any initial data $\rho_0$, we claim that there exists a $B$ such
%that $\rho_0$ obeys $\omega_B$. Indeed, according to
%\cite{kiselev2007global}, it suffices to check
%\[\omega_B\left(\frac{2\|\rho_0\|_{L^\infty}}{\|\partial_x\rho_0\|_{L^\infty}}\right)
%>2\|\rho_0\|_{L^\infty}.\]
%The inequality is clearly true if we pick
%\[B>\frac{\|\partial_x\rho_0\|_{L^\infty}}{2\|\rho_0\|_{L^\infty}}\cdot\max\left\{\delta,
%e^{\frac{2\|\rho_0\|_{L^\infty}-c}{\gamma}}\right\}.\]
%
%Then, it implies that for any $t\in[0,T]$, $\rho(\cdot,t)$ obeys
%$\omega(B)$. Hence, $\|\partial_x\rho(\cdot,t)\|_{L^\infty}\leq B$,
%and global regularity is obtained by Theorem \ref{thm:local}.
%\end{proof}
%

%The goal is to show that for any $T>0$ we may find $B>0$, so that
%$\rho_0$ obeys the rescaled
%modulus of continuity~$\omega_B$, see (\ref{oct2302}),
%and $\omega_B$
%can not be broken through before
%$T$, with the appropriate choices of~$\delta$ and $\gamma$.
Therefore, we only
need to prove Lemma~\ref{lem:oct2402}. Our strategy is as follows.
Let us assume that a modulus of continuity $\omega$, with some $\delta$ and $\gamma$
is broken at a time $t_1$, in the sense that~(\ref{oct2110}) holds
for some $x_1,y_1\in\T$. We
denote
\begin{equation}\label{oct2602}
\xi=|x_1-y_1|>0,
\end{equation}
and, for simplicity, drop the time
variable~$t_1$ in the notation. We compute:
\begin{equation}\label{eq:rhoxi}
\begin{split}
\partial_t&(\rho(x_1)-\rho(y_1))=
-\partial_x(\rho(x_1)u(x_1))+\partial_x(\rho(y_1)u(y_1))\\
&=-\big(u(x_1)\partial_x\rho(x_1)-u(y_1)\partial_x\rho(y_1)\big)
-\big(\rho(x_1)-\rho(y_1)\big)\partial_xu(x_1)
-\rho(y_1)\big(\partial_xu(x_1)-\partial_xu(y_1)\big)\\
&=I+II+III.
\end{split}
\end{equation}
We will obtain the following estimates for the three terms in the right
side of (\ref{eq:rhoxi}). To bound the first term we note that if $\Omega(\xi)$
is a modulus of continuity for $u$, then
 it follows from~\cite{kiselev2007global} that
\begin{equation}\label{oct2316bis}
|I|=|u(x_1)\partial_x\rho(x_1)-u(y_1)\partial_x\rho(y_1)|
\leq\omega'(\xi)\Omega(\xi).
\end{equation}
The modulus $\Omega(\xi)$ for $u$ is given by the following.
\begin{lem}\label{lem:oct2504}
Let $\rho$ obey the modulus of continuity $\omega$ as in (\ref{oct2102}).
There exists a universal constant~$C>0$ so that then $u(x)$ obeys a modulus of continuity
\begin{equation}\label{oct2514}
\Omega(\xi)\leq\begin{cases}C\xi,&0<\xi<\delta,\\
C\xi^{1-\alpha}{\omega(\xi)},&\xi\geq\delta.
\end{cases}
\end{equation}
\end{lem}
We will prove Lemma~\ref{lem:oct2504} later in this section.

As $\omega'(\xi)\le 1$ for $0\le\xi<\delta$, and $\omega'(\xi)=\gamma/\xi$ for $\xi>\delta$,
we conclude that
%Together, (\ref{oct2418}) and (\ref{oct2422}) show that
\begin{equation}\label{oct2424bis}
|I|\le \omega'(\xi)\Omega(\xi)\leq\begin{cases}C\xi,&0<\xi<\delta,\\
C\gamma\dfrac{\omega(\xi)}{\xi^\alpha},&\xi\geq\delta,
\end{cases}
\end{equation}
again, with the constant $C>0$ that does not depend on $\rho_0$.

%First, let~$\Omega$ be a
%modulus of continuity for $u$.
%Then, it follows from~\cite{kiselev2007global} that
%\begin{equation}\label{oct2316}
%|I|=|u(x_1)\partial_x\rho(x_1)-u(y_1)\partial_x\rho(y_1)|
%\leq\omega'(\xi)\Omega(\xi).
%\end{equation}
%Next, we will use a lower bound for $\partial_xu(x_1)=\Lambda^\alpha\rho(x_1)$
%as a function of $\xi$ only.
To bound the last two terms in the right side of (\ref{eq:rhoxi}) purely in terms of $\xi=|x_1-y_1|$ we will use the following lemma.
\begin{lem}\label{lem:oct2506} Let $\rho$ obey the modulus of continuity $\omega$ as in (\ref{oct2102}), and let $x_1$, $y_1$
be the breakthrough points as in (\ref{oct2110}).
There exists a constant $C>0$ that may only depend on $\alpha$ such that
\begin{equation}\label{oct2440bis}
\Lambda^\alpha\rho(x_1)\geq  -A(\xi),~~A(\xi):=
\begin{cases} \hbox{$C$ if $0\le\xi\le\delta$,}\\
\hbox{$ C\gamma\xi^{-\alpha}$ if $\xi>\delta$},
\end{cases}
\end{equation}
and
\begin{equation}\label{oct2322bis}
\Lambda^\alpha\rho(x_1)-\Lambda^\alpha\rho(y_1)\geq D_1(\xi),~~D_1(\xi):=\begin{cases}C\xi^{1-{\alpha}/{2}},&0<\xi\le\delta,\\
C{\omega(\xi)}{\xi^{-\alpha}},&\xi\geq\delta.
\end{cases}
\end{equation}
%with $\xi=|x_1-y_1|$.
%and
%\begin{equation}\label{oct2452bis}
%\Lambda^\alpha\rho(x_1)\geq
%- C\gamma\xi^{-\alpha}\hbox{ if $\xi>\delta$}.
%\end{equation}
\end{lem}
%
%
%\begin{equation}\label{oct2318}
%\Lambda^\alpha\rho(x_1)\geq-A(\xi).
%\end{equation}
The first estimate in the above lemma gives a bound for the second term in (\ref{eq:rhoxi}):
\begin{equation}\label{oct2320bis}
II=-\big(\rho(x_1)-\rho(y_1)\big)\Lambda^\alpha\rho(x_1)\leq\omega(\xi)A(\xi),
\end{equation}
%with $A(\xi)$ defined by the right sides of (\ref{oct2440bis}) and (\ref{oct2452bis}).
%Finally, for the third term in (\ref{eq:rhoxi}) we will use the dissipation
%bound from the following.
%\begin{lem}\label{lem:oct2510}  Let $\rho$ obey the modulus of continuity $\omega$ as in (\ref{oct2102}), and let $x_1$, $y_1$
%be the breakthrough points as in (\ref{oct2110}).
%There exists a universal constant $C>0$ so that
%\begin{equation}\label{oct2322}
%\Lambda^\alpha\rho(x_1)-\Lambda^\alpha\rho(y_1)\geq D_1(\xi):=\begin{cases}C\xi^{1-{\alpha}/{2}},&0<\xi\le\delta,\\
%C{\omega(\xi)}{\xi^{-\alpha}},&\xi\geq\delta,
%\end{cases}
%\end{equation}
%with $\xi=|x_1-y_1|$.
%\end{lem}
%This
while (\ref{oct2322bis}) leads to:
\begin{equation}\label{oct2324}
III=-\rho(y)\big(\Lambda^\alpha\rho(x_1)-\Lambda^\alpha\rho(y_1)\big)
\leq- mD_1(\xi).
\end{equation}
Here,  $m$ is the minimum of $\rho_0$ and is preserved in time: see Corollary \ref{cor:oct2502}.
%\begin{enumerate}
%\item A bound on the first term: we find a modulus of continuity for
%  $u_1$, denoted by $\Omega$. Then, it follows from
%  \cite{kiselev2007global} that
%  \[|u_1(x)\partial_x\rho(x)-u_1(y)\partial_x\rho(y)|\leq\omega'(\xi)\Omega(\xi).\]
%\item A bound on the second term: we find a lower bound of
%  $\partial_xu(x)=\Lambda^\alpha\rho(x)$
%  as a function of $\xi$, namely $\Lambda^\alpha\rho(x)\geq
%  -A(\xi)$. Then,
%  \[ -\big(\rho(x)-\rho(y)\big)\Lambda^\alpha\rho(x)\leq\omega(\xi)A(\xi).\]
%\item Dissipation on the third term: we state a dissipative estimate
%$\Lambda^\alpha\rho(x)-\Lambda^\alpha\rho(y)\geq D(\xi)$. Then,
%  \[-\rho(y)\big(\Lambda^\alpha\rho(x)-\Lambda^\alpha\rho(y)\big)\leq-\rho_mD(\xi),\]
%  where $\rho_m=\rho_m(T)$ which is defined in Lemma \ref{lem:lower}.
%\end{enumerate}
Putting (\ref{oct2316bis}),~(\ref{oct2320bis}) and~(\ref{oct2324}) together, we obtain
\begin{equation}\label{oct2326}
\partial_t(\rho(x_1,t_1)-\rho(y_1,t_1))\leq
\omega'(\xi)\Omega(\xi)+\omega(\xi)A(\xi)-mD_1(\xi).
\end{equation}
For $0\le\xi<\delta$, using (\ref{oct2514}), (\ref{oct2440bis}) and (\ref{oct2322bis}), as well as
the inequalities
\begin{equation}\label{oct2602bis}
\omega(\xi)\le \xi,~~\omega'(\xi)\le 1,~~0\le\xi<\delta,
\end{equation}
we see that
\begin{equation}\label{oct2458}
\omega'(\xi)\Omega(\xi)+\omega(\xi)A(\xi)-\frac{1}{2}mD_1(\xi)
\le C\xi-Cm\xi^{1-\alpha/2}<0,
\end{equation}
provided that
%the inequality holds if $\xi\leq
%C\rho_m\xi^{1-\frac{\alpha}{2}}$. This can be acheived if
\begin{equation}\label{oct2822}
\delta<Cm^{{2}/{\alpha}}.
\end{equation}
On the other hand,
for $\xi\geq\delta$, the above bounds tell us
\begin{equation}\label{oct2460}
\omega'(\xi)\Omega(\xi)+\omega(\xi)A(\xi)-\frac{1}{2}mD_1(\xi)
\le\farc{C\gamma\omega(\xi)}{\xi^\alpha}
-\farc{Cm\omega(\xi)}{\xi^\alpha}<0,
\end{equation}
if
\begin{equation}\label{oct2824}
\gamma<Cm.
\end{equation}
Therefore, for $\delta$ and $\gamma$ sufficiently small, we have
%\eqref{oct2454} holds
\begin{equation}\label{oct2604}
\partial_t(\rho(x_1,t_1)-\rho(y_1,t_1))<0,
\end{equation}
which is a contradiction to the assumption that $t_1$ is the first breakthrough tine.
Thus, $\omega$ can never be broken, and the proof of Lemma~\ref{lem:oct2402} is complete,
except for the proof of Lemmas~\ref{lem:oct2504} and \ref{lem:oct2506}.~$\Box$

%If we can construct a modulus $\omega$ such that
%\begin{equation}\label{eq:modulus}
%\omega'(\xi)\Omega(\xi)+\omega(\xi)A(\xi)<\frac{m}{2}D(\xi),
%\end{equation}
%then the modulus of continuity $\omega$ is preserved in time. This will be our goal.
%

\subsubsection{The dissipation bound in Lemma~\ref{lem:oct2506}}

We first prove the dissipation bound (\ref{oct2322bis}) in Lemma~\ref{lem:oct2506}.
It was shown in~\cite{kiselev2011nonlocal} that
\begin{equation}\label{oct2328}
\Lambda^\alpha\rho(x_1)-\Lambda^\alpha\rho(y_1)\geq D(\xi)
\end{equation}
with
\begin{equation}\label{oct2330}
D(\xi)=c_\alpha\left[\int_0^{{\xi}/{2}}
\frac{2\omega(\xi)-\omega(\xi+2\eta)-\omega(\xi-2\eta)}{\eta^{1+\alpha}}d\eta+
\int_{{\xi}/{2}}^\infty
\frac{2\omega(\xi)-\omega(\xi+2\eta)+\omega(2\eta-\xi)}{\eta^{1+\alpha}}d\eta
\right].
\end{equation}
Both terms in the right side are positive due to the concavity of $\omega$.

To obtain a lower bound for $D(\xi)$, we consider two cases. For
$\xi\leq\delta$, we only keep the first term. Note that
\[
\omega(\xi+ 2\eta)\le \omega(\xi)+ 2\omega'(\xi)\eta
\]
due to the concavity
of $\omega$, and
\[
\omega(\xi-2\eta)\le \omega(\xi)-2\omega'(\xi)\eta+2\omega''(\xi)\eta^2,
\]
due to the second order Taylor formula and the
monotone growth of
\[
\omega''(\xi)=-\frac{\alpha(2+\alpha)}{4}\xi^{-1+\alpha/2}.
\]
This gives
\begin{equation}\label{oct2332}
D(\xi)\geq C\int_0^{{\xi}/{2}}
\frac{(-\omega''({\xi}))\eta^2}{\eta^{1+\alpha}}d\eta=
C\xi^{1-{\alpha}/{2}},~~\hbox{ for $0\le\xi\le\delta$,}
\end{equation}
which is the first bound in (\ref{oct2322bis}).

For $\xi>\delta$, we only keep the second term in (\ref{oct2330}).
Due to the concavity of $\omega$, we
have
\begin{equation}\label{oct2340}
\omega(2\eta+\xi)-\omega(2\eta-\xi)\leq\omega(2\xi)=\omega(\xi)+\gamma\log2
\leq\frac{3}{2}\omega(\xi),
\end{equation}
if
\begin{equation}\label{oct2338}
\gamma\leq\farc{\omega(\delta)}{2\log2}=\farc{\delta-\delta^{1+\alpha/2}}
{2\log 2}.
\end{equation}
In that case, we have, using (\ref{oct2340}):
\begin{equation}\label{oct2334}
D(\xi)\geq
c_\alpha\int_{{\xi}/{2}}^\infty
\frac{2\omega(\xi)-\omega(2\xi)}{\eta^{1+\alpha}}d\eta
\geq
%c_\alpha\cdot\frac{1}{2}
C\omega(\xi)\cdot\frac{1}{\alpha}\left(\frac{\xi}{2}\right)^{-\alpha}
=C\frac{\omega(\xi)}{\xi^\alpha},~~\hbox{ for $\xi>\delta$,}
\end{equation}
and the proof of (\ref{oct2322bis}) is complete.

%Summarizing, we have shown that if (\ref{oct2338}) holds, then $D(\xi)$ obeys a
%lower bound
%\begin{equation}\label{oct2336}
%D(\xi)\geq\begin{cases}C\xi^{1-{\alpha}/{2}},&0<\xi\le\delta,\\
%C\dfrac{\omega(\xi)}{\xi^\alpha},&\xi\geq\delta.
%\end{cases}
%\end{equation}
%The constant $C$ here does not depend on the initial condition $\rho_0$.

 \subsubsection{A lower bound on $\Lambda^\alpha\rho$ in Lemma~\ref{lem:oct2506}}

The next step is obtain
the lower bound (\ref{oct2440bis}) for
$\Lambda^\alpha\rho(x_1,t_1)$.
As $\omega$ is a modulus of $\rho$, we have for any $z\in\R$
\begin{equation}\label{oct2430}
\rho(z)\leq\rho(y)+\omega(|y-z|),
\end{equation}
while
\begin{equation}\label{oct2434}
\rho(x_1)=\rho(y_1)+\omega(|x_1-y_1|).
\end{equation}
This implies a lower bound
\begin{eqnarray}\label{oct2432}
&&\Lambda^\alpha\rho(x_1)=
c_\alpha\int_\R\frac{\rho(x_1)-\rho(y_1)+\rho(y_1)-\rho(z)}
{|x_1-z|^{1+\alpha}}dz
\geq
c_\alpha\int_\R\frac{\omega(\xi)-\omega(|y_1-z|)}{|x_1-z|^{1+\alpha}}dz
\nonumber\\
&&~~~~~~~~~~~
=c_\alpha\int_\R\frac{\omega(\xi)-\omega(|\xi-\eta|)}{|\eta|^{1+\alpha}}d\eta
=:-A(\xi).
\end{eqnarray}
Our goal is to bound $A(\xi)$ from above.
%\begin{equation}\label{oct2436}
%A(\xi)\ge -C,
%\end{equation}
%with a universal constant $C>0$.
Let us decompose the integral in the second line of~(\ref{oct2432}) as
%$\displaystyle\int_\R\frac{\omega(\xi)-\omega(|\xi-\eta|)}{|\eta|^{1+\alpha}}d\eta$
%into four parts.
\[
-A(\xi)=c_\alpha\int_\R  \frac{\omega(\xi)-\omega(|\xi-\eta|)}{|\eta|^{1+\alpha}}d\eta=\int_{-\infty}^{-\xi}+\int_{-\xi}^\xi+
\int_\xi^{2\xi}+\int_{2\xi}^\infty=A_1+A_2+A_3+A_4.
\]
We claim that $A_2$ and $A_3$ are positive, so that their contribution to $A(\xi)$ is negative.
%In this four parts, \RNum{2} and \RNum{3} are postive.
Indeed, we can
estimate $A_2$ using the concavity of $\omega$:
\begin{equation}\label{oct2438}
A_2=\int_0^\xi
\frac{2\omega(\xi)-\omega(\xi-\eta)-\omega(\xi+\eta)}{\eta^{1+\alpha}}d\eta
\geq0.
\end{equation}
In addition, $A_3\geq0$ simply due to the monotonicity of $\omega$, which
implies
\[
\hbox{$\omega(\xi)\geq\omega(|\eta-\xi|)$, for $\eta\in[\xi, 2\xi]$.}
\]
It remains to bound $A_1$ and  $A_4$ from below. We first consider
$0\le\xi\leq\delta$. In this region, we can estimate~$A_4$ as follows:
\begin{eqnarray}\label{oct2440}
&&\!\!\!
A_4\geq-\int_{2\xi}^\infty\frac{\omega(\eta-\xi)}{\eta^{1+\alpha}}d\eta
\ge-\int_{2\xi}^{\xi+\delta}\frac{\eta}{\eta^{1+\alpha}}d\eta
-\int_{\xi+\delta}^\infty\frac{\gamma\log((\eta-\xi)/\delta)+
\delta-\delta^{1+\alpha/2}}
{\eta^{1+\alpha}}d\eta
\nonumber\\
&&~\geq-\int_0^{2\delta}\frac{d\eta}{\eta^{\alpha}}
-(\delta-\delta^{1+\alpha/2})\int_\delta^\infty\farc{d\eta}{\eta^{1+\alpha}}
-\gamma\int_{\delta}^\infty\frac{\log(\eta/\delta)}
{\eta^{1+\alpha}}d\eta \ge -C\delta^{1-\alpha}-C\gamma\delta^{-\alpha}.
%\nonumber\\
%&&~~~~\geq-\frac{(2\delta)^{1-\alpha}}{1-\alpha}
%-\gamma\delta^{-\alpha}\left(\frac{1}{\alpha}\log\delta+\frac{1}{\alpha^2}\right)
%-\frac{c\delta^{-\alpha}}{\alpha}.
\end{eqnarray}
%Note that
%$c\delta^{-\alpha}=\delta^{1-\alpha}-\delta^{1-\frac{\alpha}{2}}
%-\gamma\delta^{-\alpha}\log\delta$.
Thus, if we choose $\delta<1$ and
$\gamma<\delta$, as in (\ref{oct2338}),
we obtain
\begin{equation}
A_4\ge -C,~~\hbox{ for $0\le\xi\le\delta$.}
\end{equation}
The term $A_1$ can be estimated similarly for $0\le\xi\le\delta$.
Indeed, for $\xi<{\delta}/{2}$, we have
\begin{eqnarray}\label{oct2444}
&&A_1\geq-\int_\xi^\infty\frac{\omega(\eta+\xi)}{\eta^{1+\alpha}}d\eta
\ge-\int_\xi^{\delta-\xi}\frac{2\eta}{\eta^{1+\alpha}}d\eta
-\int_{\delta-\xi}^\infty\frac{\gamma\log((\eta+\xi)/\delta)+\delta}
{\eta^{1+\alpha}}d\eta\nonumber\\
&&~~~~
\geq-C{\delta^{1-\alpha}}-C
\gamma\int_{\delta/2}^\infty\frac{\log(\eta/\delta)}
{\eta^{1+\alpha}}d\eta\ge -C{\delta^{1-\alpha}}-C
\gamma\delta^{-\alpha}\ge -C,
%\nonumber\\
%&&~~~~
%\geq-\frac{\delta^{1-\alpha}}{1-\alpha}
%-2^{1+\alpha}\gamma\delta^{-\alpha}\left(\frac{1}{\alpha}\log\delta+\frac{1}{\alpha^2}\right)-\frac{2^{1+\alpha}c\delta^{-\alpha}}{\alpha}\geq-C,
\end{eqnarray}
provided that $\gamma$ satisfies (\ref{oct2338}).
On the other hand, for $\delta/2\le\xi\le\delta$, we have
\begin{eqnarray}\label{oct2446}
&&A_1\geq
-\int_\xi^\infty\frac{\gamma\log((\eta+\xi)/\delta)+\delta}{\eta^{1+\alpha}}d\eta
\geq -\int_\xi^\infty\frac{\gamma\log(2\eta/\delta)+\delta}{\eta^{1+\alpha}}d\eta
\geq -\int_{{\delta}/{2}}^\infty
\frac{\gamma\log(2\eta/\delta)+\delta}{\eta^{1+\alpha}}d\eta\nonumber
\\
&&~~~~\ge -C\delta^{1-\alpha}-C\gamma\delta^{-\alpha}\ge -C.
%&&~~~~\geq
%-\gamma\left(\frac{\delta}{2}\right)^{-\alpha}
%\left[\frac{1}{\alpha}\log\left(\frac{\delta}{2}\right)
%+\frac{1}{\alpha^2}\right]-\frac{2^\alpha c\delta^{-\alpha}}{\alpha}\geq-C.
\end{eqnarray}
%For the third inequality, we use the fact that $\gamma\log(2\eta)+c>0$
%for $\eta\geq\frac{\delta}{2}$, provided that $\delta<\frac{1}{2}$ and
%$\gamma<\frac{1-2^{-\frac{\alpha}{2}}}{2\log2}\delta$.
%
Summing up the above computation, we conclude that
\begin{equation}\label{oct2440bis3}
\hbox{$\Lambda^\alpha\rho(x_1)\geq -A(\xi)\ge -C$ if $0\le\xi\le\delta$.}
\end{equation}
%As $\omega(\xi)\leq\xi$, we deduce
%\begin{equation}\label{oct2442}
%\hbox{$\omega(\xi)A(\xi)\geq -C\xi$
%for $0\le\xi\le\delta$.}
%\end{equation}
%
On the other hand, if $\xi>\delta$,
we have the following estimates on $A_1$ and $A_4$:
\begin{equation}\label{oct2448}
A_{1}=\gamma\int_{-\infty}^{-\xi}
\frac{\log\xi-\log(\xi-\eta)}{|\eta|^{1+\alpha}}d\eta
=-\frac{\gamma}{\xi^\alpha}
\int_{-\infty}^{-1}\frac{\log(1-\zeta)}{|\zeta|^{1+\alpha}}d\zeta\geq
-\frac{C\gamma}{\xi^\alpha},
\end{equation}
and
\begin{equation}\label{oct2450}
A_4=\gamma\int_{2\xi}^\infty
\frac{\log\xi-\log(\eta-\xi)}{|\eta|^{1+\alpha}}d\eta
=-\frac{\gamma}{\xi^\alpha}
\int_2^\infty\frac{\log(\zeta-1)}{\zeta^{1+\alpha}}d\zeta\geq
-\frac{C\gamma}{\xi^\alpha},
\end{equation}
Thus, we have the bound
\begin{equation}\label{oct2452}
\Lambda^\alpha\rho(x_1)\geq -A(\xi)\ge
- C\gamma\xi^{-\alpha}\hbox{ for $\xi>\delta$},
\end{equation}
finishing the proof of (\ref{oct2440bis}), as well as of Lemma~\ref{lem:oct2506}.~$\Box$

%Put everything together, we conclude with
%\[
%\omega(\xi)A(\xi)\leq\begin{cases}C\xi&\xi<\delta,\\
%C\gamma\frac{\omega(\xi)}{\xi^\alpha}&\xi\geq\delta.
%\end{cases}
%\]

\subsubsection{The proof of Lemma~\ref{lem:oct2504}}

Next, we find a modulus of continuity $\Omega$ for $u$, if $\rho$ obeys $\omega$ given by (\ref{oct2102}).
We start with~(\ref{oct2342}):
\begin{equation}\label{oct2402}
u(x) %=c_\alpha p.v. \int\frac{\varphi(x)-\varphi(x+y)}{|y|^{1+\alpha}}\, dy
=c_\alpha \lim_{\epsilon\downarrow 0} \int_{|y|>\epsilon}
\frac{\varphi(x)-\varphi(x+y)}{|y|^{1+\alpha}}\, dy.
\end{equation}
The first term in the right side can evaluated explicitly:
\begin{equation}\label{oct2404}
\int_{|y|>\epsilon}  \frac{\varphi(x)}{|y|^{1+\alpha}}\, dy
%=\frac{2}{\alpha}\frac{\varphi(x)}{\epsilon^{\alpha}} +\frac{1}{\alpha}
%\int_{|y|>\epsilon} \frac{\theta(x)}{sgn(y)|y|^\alpha}\, dy
=\frac{2}{\alpha}\frac{\varphi(x)}{\epsilon^{\alpha}}.
\end{equation}
%where the second equality is due to symmetry.
The second term in the right side of (\ref{oct2402}) can be re-written
using integration by parts as
\begin{eqnarray}\label{oct2406}
\int_{|y|>\epsilon} \frac{\varphi(x+y)}{|y|^{1+\alpha}}\, dy = \frac{1}{\alpha} \frac{\varphi(x+\epsilon)+\varphi(x-\epsilon)}{\epsilon^\alpha} + \frac{1}{\alpha} \int_{|y|>\epsilon} \frac{\theta(x+y)}{\hbox{sgn}(y)|y|^\alpha}\, dy.
\end{eqnarray}
As $\theta\in L^\infty$, so that $\varphi$ is uniformly Lipschitz,
we can combine (\ref{oct2404}) and (\ref{oct2406}), pass to the
limit~$\eps\downarrow 0$, and obtain
\begin{equation}\label{oct2408}
u(x) = -\frac{c_\alpha}{\alpha}\int_\R
\frac{\theta(x+y)}{\hbox{sgn}(y)|y|^\alpha}\, dy.
\end{equation}
Let us note that, since $\theta(x)$ is a periodic mean-zero function, the
integral in the right side of~(\ref{oct2408}) converges as $|y|\to+\infty$,
and
% \begin{eqnarray}\label{oct2410}
% &&u(x) = -\lim_{\epsilon\downarrow 0}
% \frac{c_\alpha}{\alpha}\int_{|y|\le\epsilon^{-1}}
% \frac{\theta(x+y)}{\hbox{sgn}(y)|y|^\alpha}\, dy=
% \lim_{\epsilon\downarrow 0}
% \frac{c_\alpha}{\alpha}\int_{|y|\le\epsilon^{-1}}
% \frac{\theta(x)-\theta(x+y)}{\hbox{sgn}(y)|y|^\alpha}\, dy\nonumber\\
% &&~~~~~~=\lim_{\epsilon\downarrow 0}
% \frac{c_\alpha}{\alpha}\int_{|y|\le\epsilon^{-1}}
% \frac{\rho(x)-\rho(x+y)}{\hbox{sgn}(y)|y|^\alpha}\, dy.
% \end{eqnarray}
\begin{equation}\label{oct2410}
u(x) =  \frac{c_\alpha}{\alpha}\int_\R
\frac{\theta(x)-\theta(x+y)}{\hbox{sgn}(y)|y|^\alpha}\, dy
=\frac{c_\alpha}{\alpha}\int_\R  \frac{\rho(x)-\rho(x+y)}{\hbox{sgn}(y)|y|^\alpha}\, dy.
\end{equation}
%Observing that $\int_\R sgn(y)|y|^{-\alpha}dy=0$, we then express $u$
%in the following form in terms of $\rho$:
%\begin{align}
%\label{u}
%u(x) = \frac{c_\alpha}{\alpha}\int_R \frac{\rho(x)-\rho(x+y)}{sgn(y)|y|^\alpha}\, dy.
%\end{align}
Using an argument similar to that in the appendix
of \cite{kiselev2007global}, one can show that, as long as~$\rho(x)$ obeys a modulus of
continuity $\omega$, the function
$u(x)$ given by~\eqref{oct2410} obeys the modulus of continuity
\begin{equation}\label{oct2412}
\Omega(\xi) = C\left(\int_0^\xi \frac{\omega(\eta)}{\eta^\alpha}\, d\eta
+ \xi \int_\xi^\infty \frac{\omega(\eta)}{\eta^{1+\alpha}}\,
d\eta \right),
\end{equation}
with a universal constant $C>0$.

Thus, for $0\le\xi\leq\delta$, we get
\begin{align}\label{oct2612}
\Omega(\xi)\leq& C\left(\int_0^\xi\eta^{1-\alpha}d\eta
+\xi\int_\xi^\delta\eta^{1-\alpha}d\eta
+\xi\int_\delta^\infty\frac{\gamma\log(\eta/\delta)+\delta }{\eta^{1+\alpha}}d\eta
\right)\nonumber\\
\leq& C\left(\xi^{2-\alpha}+\xi\delta^{2-\alpha}+\xi\gamma\delta^{-\alpha}+
\xi\delta^{1-\alpha}\right)\le C\xi, %(\gamma|\log\delta|+c)+\xi\gamma\right).
\end{align}
as long as we take $\gamma<\delta$.
%with some $p>1$, then (\ref{oct2416})
%implies that $|c|\le C\delta$, thus
%\[
%\delta^{-\alpha}(|\gamma\log\delta|+|c|)\le C.
%\]
%Therefore,
%we obtain
%\begin{equation}\label{oct2610}
%\hbox{$\Omega(\xi)\leq C\xi$ for $0\le\xi\le\delta$,}
%\end{equation}
This is the first inequality in (\ref{oct2514}).
%
%and since $\omega'(\xi)\leq1$, we conclude that
%\begin{equation}\label{oct2418}
%\omega'(\xi)\Omega(\xi)\leq C\xi.
%\end{equation}

For $\xi>\delta$, we use % (\ref{oct2416}) and
(\ref{oct2412}) to write
\begin{eqnarray}\label{oct2420}
&&\Omega(\xi)\leq C\Big(\int_0^\delta\eta^{1-\alpha}d\eta
+\int_\delta^\xi\frac{\gamma\log(\eta/\delta)
+\delta-\delta^{1+\alpha/2}}{\eta^\alpha}d\eta
+\xi\int_\xi^\infty\frac{\gamma\log(\eta/\delta)
+\delta-\delta^{1+\alpha/2}}{\eta^{1+\alpha}}d\eta
\Big)\nonumber\\
&&~~~~~~
\le C\left( \delta^{2-\alpha}+\xi^{1-\alpha}(\delta-\delta^{1+\alpha/2})
\right)+C\gamma\delta^{1-\alpha}\int_1^{\xi/\delta}\frac{\log\eta}{\eta^\alpha}
d\eta+C\gamma\xi\delta^{-\alpha}\int_{\xi/\delta}^\infty\farc{\log\eta}
{\eta^{1+\alpha}}d\eta\\
&&~~~~~~
\le C\left( \delta^{2-\alpha}+\xi^{1-\alpha}(\delta-\delta^{1+\alpha/2})
\right)+C\gamma\xi^{1-\alpha}(1+\log(\xi/\delta)) \leq C\left(\delta^{2-\alpha}
+\xi^{1-\alpha}\omega(\xi)\right)\le C\xi^{1-\alpha}\omega(\xi),\nonumber
%\nonumber\\
%&&\leq C\left(\delta^{2-\alpha}+\xi^{1-\alpha}(\gamma\log(\xi/\delta)+
%\delta-\delta^{1+\alpha/2})\right)
%\leq C\left(1+\xi^{1-\alpha}\omega(\xi)\right).
\end{eqnarray}
finishing the proof of Lemma~\ref{lem:oct2504}.

\subsection{The global regularity for the full system}

We now consider the full system (\ref{eq:rhobis})-(\ref{eq:uxbis})
\begin{eqnarray}
&&\partial_t\rho+\partial_x(\rho u)=0,\label{eq:rhobis2}\\
&&\partial_tG+\partial_x(Gu)=0,\label{eq:Gbis2}\\
 &&\partial_xu=\Lambda^\alpha\rho+G,\label{eq:uxbis2}
\end{eqnarray}
without the extra assumption $G\equiv 0$.
Let us recall representation (\ref{oct1714}):
\begin{equation}\label{oct2202}
u(x)=\Lambda^\alpha\varphi(x)+(\psi(x)+I_0)=:u^{(1)}(x)+u^{(2)}(x).
\end{equation}
Here, $\phi(x)$ and $\psi(x)$ are the mean-zero primitives of $\theta$ and $G$, respectively, as in
(\ref{oct1724})-(\ref{oct1726}), and~$I_0$ is given by (\ref{eq:C0}).
%The contributions of  $u_1$ and $u_2$ will be discussed separately.
%The strategy of the global regularity proof is as in the case $G\equiv 0$, except for an additional
%scaling argument.

Note that if $\rho(x,t)$ and $G(x,t)$ are solutions of (\ref{eq:rhobis2})-(\ref{eq:uxbis2}),
with the corresponding velocity~$u(x,t)$, then
\begin{equation}\label{oct2640}
\rho_\lambda(x,t)=\rho(\lambda x,\lambda^\alpha t),~~G_\lambda(x,t)=\lambda^\alpha G(\lambda x,\lambda^\alpha t),
\end{equation}
are also solutions, with the corresponding velocity
\begin{equation}\label{oct2642}
u_\lambda(x,t)=\lambda^{-(1-\alpha)}u(\lambda x,\lambda^\alpha t),
\end{equation}
and
\begin{equation}\label{oct2644}
F_\lambda(x,t)=\lambda^\alpha F(\lambda x,\lambda^\alpha t),~~F(x,t)=\frac{G(x,t)}{\rho(x,t)}.
\end{equation}
Note that if $\rho_\la(x,t)$ obeys a modulus of continuity $\omega$, then $\rho(x,t)$ obeys the modulus
of continuity
\begin{equation}\label{oct2802}
\omega_\lambda(\xi)=\omega(\lambda^{-1}\xi).
\end{equation}
The proof of the global regularity for the solutions of (\ref{eq:rhobis2})-(\ref{eq:uxbis2}) is based on the following lemma.
\begin{lem}\label{lem:oct2802}
Let $\omega$ and $\omega_\lambda$ be as in (\ref{oct2102}) and (\ref{oct2802}), respectively.
Given a smooth periodic initial condition $(\rho_0,u_0)$ for (\ref{eq:rhobis2})-(\ref{eq:uxbis2}) and $T>0$, there
exist $\delta>0$, $\gamma>0$ and $\lambda>0$ so that $\rho(x,t)$ obeys the modulus of continuity $\omega_\lambda(\xi)$ for all
$0\le t\le T$. The parameters $\delta$, $\gamma$ and~$\la$ may depend on~$\alpha,$ $\rho_0,$ $u_0,$ and $T$.
\end{lem}
This will imply a uniform bound on $\|\partial_x\rho\|_{L^\infty}$ on $0\le t\le T$. As $T$ is arbitrary,
this is sufficient for the global regularity of the solutions, according to~\eqref{eq:BKM}.
Note that   $\rho(x,t)$ obeys $\omega_\lambda$ until a time $T$ if and only if $\rho_\lambda(x,t)$
obeys the modulus of continuity~$\omega$ until the time~$T_\lambda=\lambda^{-\alpha}T$, and this is what we will show.
That is, given $\rho_0$ and $u_0$, and $T>0$, we will find $\lambda>0$, $\delta>0$ and $\gamma>0$ sufficiently small, so that
(i) $\rho_\lambda(0,x)=\rho_0(\lambda x)$ obeys~$\omega$, and~(ii)~$\rho_\la(x,t)$ obeys $\omega$ at least until the time $\lambda^{-\alpha}T$.
The a priori bounds on $\rho(x,t)$ and~$F(x,t)$ will play a crucial role in the proof.

As in the case $G\equiv 0$ considered above, we assume that a modulus of continuity~$\omega$ of the form~ (\ref{oct2102}),
with some $\delta$ and $\gamma$,
is broken by $\rho_\lambda$ at a time $t_1$, at some $x_1,y_1\in\Rm$, in the sense of~(\ref{oct2110}). If $\T=[0,L]$,
then $\rho_\lambda$ is $\lambda^{-1}L$-periodic,
and we can restrict our attention to~$x_1,y_1\in\T_\lambda:=\lambda^{-1}\T$. We also   set
\begin{equation}\label{oct2620}
\xi=|x_1-y_1|>0,
\end{equation}
and drop the time
variable~$t_1$ in the notation. We decompose as in (\ref{eq:rhoxi}):
\begin{equation}\label{oct2622}
\partial_t(\rho_\lambda(x_1)-\rho_\lambda(y_1))=
-\partial_x(\rho_\la(x_1)u_\la(x_1))+\partial_x(\rho_\la(y_1)u_\la(y_1))=R_1+R_2,
\end{equation}
with the terms $R_1$ and $R_2$ coming from the contributions of $u_\lambda^{(1)}$ and $u_\lambda^{(2)}$ in (\ref{oct2202}).
We treat~$R_1$ as before:
\begin{equation}
\begin{split}
R_1&=-\big(u_\la^{(1)}(x_1)\partial_x\rho_\la(x_1)-u_\la^{(1)}(y_1)\partial_x\rho_\la(y_1)\big)\\
&-\big(\rho_\la(x_1)-\rho_\la(y_1)\big)\partial_xu_\la^{(1)}(x_1)
-\rho_\la(y_1)\big(\partial_xu_\la^{(1)}(x_1)-\partial_xu_\la^{(1)}(y_1)\big)
=I+II+III.
\end{split}
\end{equation}
%\begin{equation}\label{oct2622}
%\begin{split}
%\partial_t&(\rho(x_1)-\rho(y_1))=
%-\partial_x(\rho(x)u(x_1))+\partial_x(\rho(y_1)u(y_1))\\
%&=-\big(u(x_1)\partial_x\rho(x_1)-u(y_1)\partial_x\rho(y_1)\big)
%-\big(\rho(x_1)-\rho(y_1)\big)\partial_xu(x_1)
%-\rho(y_1)\big(\partial_xu(x_1)-\partial_xu(y_1)\big)\\
%&=I_1+I_2+II_1+II_2+III_1+III_2.
%\end{split}
%\end{equation}
%We have split in the last equality the contributions of $u_1$ and $u_2$ to $I$, $II$ and $III$.
Note that $I$ and~$II$ can be estimated exactly as before: first, as in (\ref{oct2424bis}), we have
\begin{equation}\label{oct2628}
|I|\le  \begin{cases}C\xi,&0<\xi<\delta,\\
C\gamma\dfrac{\omega(\xi)}{\xi^\alpha},&\xi\geq\delta,
\end{cases}
\end{equation}
with a constant $C>0$ that does not depend on $\rho_0$ or $u_0$. The term $II$ can be bounded
as in~(\ref{oct2320bis}):
\begin{equation}\label{oct2630bis}
II \leq\omega(\xi)A(\xi),
\end{equation}
with $A(\xi)$ defined in (\ref{oct2440bis}). The term $III$ is bounded slightly differently from (\ref{oct2324})
\begin{equation}\label{oct2632}
III\leq- \rho_m^{(\lambda)}(T)D_1(\xi).
\end{equation}
Here, $\rho_m^{(\lambda)}(T)$ is the minimum of $\rho_\lambda(x,t)$ over $0\le t\le \lambda^{-\alpha}T$,
and $D_1(\xi)$  is defined in (\ref{oct2322bis}).  The lower bound (\ref{oct2810}) in
Lemma~\ref{lem:lower} implies that
\begin{equation}\label{oct2812}
\begin{split}
\rho_m^{(\lambda)}(T)\ge&\farc{1}{[\rho_m^{(\lambda)}(0)]^{-1}+\lambda^{-\alpha}T\|F_0^{\lambda}\|_{L^\infty}}=
\farc{1}{[\rho_m(0)]^{-1}+ T\|F_0\|_{L^\infty}}\\
\ge &\farc{\rho_m(0)}{1+T\|\partial_xu_0\|_{L^\infty}+T\|\Lambda^\alpha\rho_0\|_{L^\infty}}:=\bar\rho_m(T),
\end{split}
\end{equation}
as follows from (\ref{oct2644}).
That is, even though now, unlike in the special case $G\equiv 0$, the function~$\rho(x,t)$
does not necessarily obey the minimum principle, and $\rho_m(t)$ may decrease in time, the value of~$\rho_m^{(\lambda)}(t)$
does not depend on $\lambda>0$. Thus, we may first choose the parameters~$\delta$ and $\gamma$ in the definition~(\ref{oct2102})
of the modulus of continuity $\omega$ so that (\ref{oct2822}) and~(\ref{oct2824}) hold with $m$ replaced by~$\bar\rho_m(T)$, and,
in addition, they satisfy  (\ref{oct2338}). Next, we choose $\lambda$ sufficiently small, so that $\rho_\lambda^0(x)=\rho_0(\lambda x)$
obeys the modulus of continuity $\omega$ with the above choice of $\delta$ and $\gamma$.

It remains to take into account the contribution
of $u_\lambda^{(2)}$ to the right side of~(\ref{oct2622}).
%
%
%
%Now, we consider the full system \eqref{eq:rhoG}, taking into account
%of the contribution from $u_2$.
The goal is to control the corresponding  terms in \eqref{eq:rhoxi} by the
dissipation, namely, to show that
\begin{eqnarray}\label{oct2826}
&&R_2=\big|u_\lambda^{(2)}(x_1)\partial_x\rho_\lambda(x_1)-u_\lambda^{(2)}(y_1)\partial_x\rho_\lambda(y_1)\big|
+\big|\rho_\la(x_1)\partial_xu_\la^{(2)}(x_1)-\rho_\la(y_1) \partial_xu_\la^{(2)}(y_1)\big|\nonumber\\
&&~~~~=R_{21}+R_{22}
 <\frac{1}{2}\bar\rho_m(T)D_1(\xi).
\end{eqnarray}
%We directly work with $\omega_B$. The corresponding dissipative
%estimate reads
%\[D(\xi)\geq\begin{cases}
%CB^{1+\frac{\alpha}{2}}\xi^{1-\frac{a}{2}}&\xi<\frac{\delta}{B},\\
%C\frac{\omega_B(\xi)}{\xi^\alpha}&\frac{\delta}{B}\leq\xi\leq d,
%\end{cases}\]
%where $d=|\T|$ is the size of the period. We do not need to discuss
%$\xi>d$ as the breakthrough can not happen at such $\xi$.
Note that the flow $u_\la^{(2)}(x)$ is Lipschitz,
as
\begin{equation}\label{oct2820}
|\partial_xu_\la^{(2)}(t,x)|=|G_\la(t,x)|\le|\rho_\lambda(t,\cdot)\|_{L^\infty}\|F_\lambda(t,\cdot)\|_{L^\infty}\le C_0\lambda^\alpha,
\end{equation}
with a constant $C_0$ that depends on the initial conditions $\rho_0$ and $u_0$ but not on $\lambda>0$.
Therefore,~$u_\lambda^{(2)}$ obeys the modulus of continuity
\begin{equation}\label{oct2828}
\Omega_2(\xi)=C_0\lambda^\alpha\xi,%\|G\|_{L^\infty}\xi.
\end{equation}
%It implies
and the first term in (\ref{oct2826}) can be bounded by
\begin{equation}\label{oct2830}
R_{21}:=\big|u_\la^{(2)}(x_1)\partial_x\rho_\la(x_1)-u_\la^{(2)}(y_1)\partial_x\rho_\lambda(y_1)\big|\leq
C_0\lambda^\alpha\xi\omega'(\xi). %\leq\frac{1}{4}\rho_mD(\xi).
\end{equation}
Let us recall from (\ref{oct2102}) and (\ref{oct2322bis}) that
\begin{equation}\label{oct2834}
\omega'(\xi)\le 1,~~D_1(\xi)=C_1\xi^{1-\alpha/2},~~\hbox{ for $0\le\xi\le\delta$},
\end{equation}
hence, we have
%choose $\delta$ small enough and $B>1$,
\begin{align}\label{oct2832}
R_{21}\le C_0\lambda^\alpha\xi \omega'(\xi) \leq C_0 \lambda^\alpha \xi < \frac{C_1\bar\rho_m(T)}{4}
\xi^{1-\frac{\alpha}{2}} < \frac{1}{4} \bar\rho_m(T) D_1(\xi) ~~\hbox{ for $0\le\xi\le\delta$},
\end{align}
provided that $\delta$ and $\lambda$ are sufficiently small.
%\begin{align*}
%C_0\xi \omega_B'(\xi) \leq C_1 B \xi < \frac{C}{4}
%  B^{1+\frac{\alpha}{2}}
%\xi^{1-\frac{\alpha}{2}} < \frac{1}{4} \rho_m D(\xi). ~~\hbox{ for $0\le\xi\le\delta$},
%\end{align*}
On the other hand,  we see from (\ref{oct2102}) and (\ref{oct2322bis}) again that
\begin{equation}\label{oct3012}
\omega'(\xi)=\farc{\gamma}{\xi},~~D_1(\xi)=\frac{C_1\omega(\xi)}{\xi^\alpha},~~\hbox{ for $\delta\le\xi\le L\lambda^{-1}$}.
\end{equation}
%in the region $\delta\leq\xi\leq  $,
%using $\omega'(\xi)= \frac{\gamma}{\xi}$ and
It is also straightforward to check that $D_1(\xi)$ is decreasing for $\xi>\delta$, provided that
\begin{equation}\label{oct2834bis}
\gamma<c\delta,
\end{equation}
with a sufficiently small constant $c>0$ that depends only on $\alpha$.
We also have
\begin{equation}\label{oct2838}
\frac{\omega(\lambda^{-1}L)}{L^\alpha}\to +\infty,~~\hbox{as $\lambda\to 0$, with $L>0$ fixed}.
\end{equation}
Hence, %possibly decreasing $\gamma$ further, and
taking $\lambda$ sufficiently small, depending on $L$ as well,
we have the inequality
\begin{equation}\label{oct2836}
R_{21}\le C_0\lambda^\alpha\xi\omega'(\xi) = C_0\lambda^\alpha\gamma<
\frac{C_1\bar\rho_m(T)}{4}\frac{\omega(\lambda^{-1}L)}{(\lambda^{-1}L)^\alpha}\leq
\frac{C_1\bar\rho_m(T)}{4}\frac{\omega(\xi)}{\xi^\alpha}
\leq\frac{1}{4}\bar\rho_m(T)D_1(\xi),
\end{equation}
for $\delta\le\xi\le L\lambda^{-1}$. Together, (\ref{oct2832}) and (\ref{oct2836}) show that
\begin{equation}\label{oct3020}
R_{21}\leq\frac{1}{4}\bar\rho_m(T)D_1(\xi).
\end{equation}

For the second term in (\ref{oct2826}), we write
\begin{align}\label{oct2840}
R_{22}&=\big|\rho_\la(x_1)\partial_xu_\la^{(2)}(x_1)-\rho_\la(y_1) \partial_xu_\la^{(2)}(y_1)\big|=
\big|\rho_\la(x_1)^2F_\la(x_1)-\rho_\la(y_1)^2F_\la(y_1)\big|\nonumber\\
&\leq
2\la^\alpha\|\rho\|_{L^\infty}^2\|F\|_{L^\infty}\leq C_0\la^\alpha,
\end{align}
with a constant $C_0$  that depends only on the initial condition $\rho_0$ and $u_0$.
Then, for $\lambda$ sufficiently small,  we have, once again using the fact that $\omega(\xi)/\xi^\alpha$
is decreasing for $\xi>\delta$ and (\ref{oct2838}):
\begin{equation}\label{oct3002}
%\big|\rho(x)\partial_xu_2(x)-\rho(y) \partial_xu_2(y)\big| \leq C_2<
R_{22}\le C_0\lambda^\alpha\le \frac{C_1\bar\rho_m(T)}{4}\frac{\omega(\lambda^{-1}L)}{(\lambda^{-1}L)^\alpha}\leq
\frac{C_1\bar\rho_m(T)}{4}\frac{\omega(\xi)}{\xi^\alpha}
=\frac{1}{4}\bar\rho_m(T)D_1(\xi),\hbox{  for $\delta\leq\xi\leq \lambda^{-1}L$. }
\end{equation}
%where the second inequality can be obtain by picking $B$ large enough.
To bound $R_{22}$ in the region $0\le\xi\le\delta$, we write
\begin{align}\label{oct3004}
R_{22}&=\big|\rho_\la(x_1)^2F_\la(x_1)-\rho_\la(y_1)^2F_\la(y_1)\big|
\\
&\le
\big|\rho_\la(x_1)^2F_\la(x_1)-\rho_\la(y_1)^2F_\la(x_1)\big|+\big|\rho_\la(y_1)^2F_\la(x_1)-\rho_\la(y_1)^2F_\la(y_1)\big|
\nonumber\\
&\le 2\|\rho_\la\|_{L^\infty}\|F_\la\|_{L^\infty}\omega(\xi)+\|\rho_\la\|_{L^\infty}^2\|\partial_xF_\la\|_{L^\infty}\xi
\nonumber
\end{align}
Lemma~\ref{lem:oct2604} guarantees that
$F$ is Lipschitz, and the Lipschitz bound is uniform in time, thus~(\ref{oct2644}) implies
%\begin{equation}\label{oct3008}
\[ \|\partial_xF_\la\|_{L^\infty}\le C_0\lambda^{1+\alpha},\]
%\end{equation}
with a constant $C_0$ that depends only on the initial conditions. In addition, it follows from~(\ref{oct2644})
that
%\begin{equation}\label{oct3010}
\[ \|F_\la\|_{L^\infty}\le C_0\lambda^\alpha. \]
%\end{equation}
Inserting the last two bounds in (\ref{oct3004}), together with the expression for $D_1(\xi)$ in (\ref{oct2834}), gives
\begin{equation}\label{oct3006}
R_{22}\le C_0\lambda^\alpha(\omega(\xi)+\xi)\le\frac{C_1\bar\rho_m(T)}{4}\xi^{1-\alpha/2}=\farc{\bar\rho_m(T)}{4}D_1(\xi).
\end{equation}
Here the constant $C_0$ depends only on the initial conditions $\rho_0$ and $u_0$, and the second inequality holds provided that
$\delta$ and $\lambda$ are sufficiently small. This proves (\ref{oct2826}), and finishes the proof of Lemma~\ref{lem:oct2802}.

Let us recap the order in which we choose the parameters.
The value of $\alpha$ is fixed throughout the argument. Given the initial data, we also fix its period, $L.$
We can also assume that~$\lambda$ does not exceed one.
Next we choose $\delta$ sufficiently small so that \eqref{oct2822} (with $m$ replaced
by~$\bar \rho_m(T)$), \eqref{oct2832}, and \eqref{oct3006} hold.
Then we choose $\gamma$ so that \eqref{oct2824} (with $m$ replaced by~$\bar \rho_m(T)$), \eqref{oct2338} and \eqref{oct2834bis} hold.
Finally, we choose $\lambda$
so that $\rho_\lambda (0,x)$ obeys $\omega$ with the above choice of $\delta,$ $\gamma$ and so that \eqref{oct2836} and \eqref{oct3002} hold.
%Once again, $\delta>0$ can be chosen first, independently of the choice of $\lambda>0$.
The proof of Theorem~\ref{thm:global}
is now complete.~$\Box$

%\begin{lem}
%$F$ is Lipschitz, and the Lipschitz bound is uniformly in time.
%\end{lem}
%\begin{proof}
%Recall that $F$ satisfies
%\[\partial_tF+u\partial_xF=0.\]
%Apply $\partial_x$ to the equation and get
%\[\partial_t\partial_xF+\partial_x(u\partial_xF)=0.\]
%So $\partial_xF$ satisfies the same continuity equation as $\rho$.
%Let $w=F/\rho$. Then, we have
%\[\partial_tw+u\partial_xw=0.\]
%It implies $\|w(\cdot,t)\|_{L^\infty}=\|w_0\|_{L^\infty}$, and
%therefore,
%\[\|\partial_xF(\cdot,t)\|_{L^\infty}\leq \|w_0\|_{L^\infty}\|\rho(\cdot,t)\|_{L^\infty}.\]
%It then follows from Theorem \ref{thm:linf} that $F$ is Lipschitz with
%a time-independent Lipschitz bound.
%\end{proof}
%
%Applying the lemma, we have
%\begin{align*}
%\big|\rho(x)\partial_xu_2(x)&-\rho(y) \partial_xu_2(y)\big|=
%\big|\rho(x)^2F(x)-\rho(y)^2F(y)\big|\\
%\leq&2\|\rho\|_{L^\infty}\|F\|_{L^\infty}\omega_B(\xi)
%+\|\rho\|_{L^\infty}^2\|\partial_xF\|_{L^\infty}\xi
%\leq C_3(\omega_B(\xi)+\xi).
%\end{align*}
%Choosing $\delta$ small enough and $B$ large enough, we conclude
%\begin{align*}
%C_3(\omega_B(\xi)+\xi) < C_3(B+1)\xi <
%\frac{C}{4} B^{1+\frac{\alpha}{2}}\xi^{1-\frac{\alpha}{2}} <\frac{1}{4} \rho_m D(\xi).
%\end{align*}
%
%Knowing that the contribution from $u_2$ is controled by dissipation,
%we use the same argument in the proof of Theorem \ref{thm:special} to
%conclude the global wellposedness of \eqref{eq:main}.

\appendix
\section{The proof of a commutator estimate}\label{sec:app}

In this section, we prove the commutator estimate \eqref{est:alpha},
\[
\|[\Lambda^\gamma, f,g]\|_{L^2}\lesssim
\|f\|_{L^2}\|g\|_{C^{\gamma+\epsilon}},\quad \gamma\in(0,1).
\]
The proof is for $x\in\R^n$, though it can be easily adapted to periodic case.
%We follow the notations in the book
%\cite{bahouri2011fourier}, summarized below.
%
Let $(\chi, \eta)$ be smooth functions such that $\chi$ is supported in
a ball $\{\xi~:~|\xi|\leq4/3\}$, $\eta$ is supported in an annulus
 $\{\xi~:~3/4\leq|\xi|\leq8/3\}$, and
\[\chi(\xi)+\sum_{q=0}^\infty\eta(2^{-q}\xi)\equiv1,\quad\forall~\xi\in\R^n.\]
It is standard to take
\[ \eta(\xi) = \chi(\xi/2) - \chi(\xi), \]
which we will assume.
Denote the Littlewood-Paley decomposition of $f$ as
$\sum_{q=-1}^\infty\Delta_qf$, where
$\Delta_qf=(\eta(2^{-q}\xi)\hat{f}(\xi))^\vee$ for $q\geq0$,
and $\Delta_{-1}f=(\chi(\xi)\hat{f}(\xi))^\vee$.
The Besov norm is defined as~\cite{bahouri2011fourier}
\[\|f\|_{B^s_{p,r}}=\left(\sum_q2^{rs}\|\Delta_qf\|_{L^p}^r\right)^{1/r}.\]

Let the partial sum
$S_qf=\sum_{p\leq q-1}\Delta_pf$. The Bony
decomposition states
\[fg=T_fg+T_gf+R(f,g),\]
where
\[T_fg=\sum_qS_{q-1}f\cdot\Delta_qg,
\quad R(f,g)=\sum_q\tilde{\Delta}_qf\cdot\Delta_qg,
\quad\tilde{\Delta}_qf=\sum_{p=q-1}^{q+1}\Delta_pf.\]

\begin{proof}[Proof of the commutator estimate]
First, we observe
\[\|f\Lambda^\gamma g\|_{L^2}\leq\|f\|_{L^2}\|\Lambda^\gamma
g\|_{L^\infty}\lesssim\|f\|_{L^2}\|g\|_{C^{\gamma+\epsilon}}.\]
Therefore, it suffies to prove
\[\|\Lambda^\gamma(fg)-g\Lambda^\gamma f\|_{L^2}\lesssim
\|f\|_{L^2}\|g\|_{C^{\gamma+\epsilon}}.\]
We apply the Bony decomposition to both terms, to get
\begin{align*}
\Lambda^\gamma(fg)=&\Lambda^\gamma(T_fg)+\Lambda^\gamma(T_gf)+\Lambda^\gamma(R(f,g))
=\RNum{1}_1+\RNum{1}_2+\RNum{1}_3,\\
g\Lambda^\gamma f=&T_{(\Lambda^\gamma f)}g+
T_g(\Lambda^\gamma f)+R(\Lambda^\gamma f,g)
=\RNum{2}_1+\RNum{2}_2+\RNum{2}_3.
\end{align*}
The terms
$\RNum{1}_1,\RNum{2}_1,\RNum{1}_3,\RNum{2}_3$ can be estimated with
standard paraproduct calculus, sketched as follows.
\begin{align*}
\|\RNum{1}_1\|_{L^2}^2=&\sum_q\|\Delta_q\Lambda^\gamma
  (T_fg)\|_{L^2}^2\lesssim\sum_q 2^{2q\gamma}\|\Delta_q(T_fg)\|_{L^2}^2
\lesssim\sum_q 2^{2q\gamma}\|S_{q-1}f\cdot\Delta_qg\|_{L^2}^2
\leq\|f\|_{L^2}^2\|g\|_{B^\gamma_{\infty,2}}^2,\\
\|\RNum{2}_1\|_{L^2}^2=&
\sum_q\|\Delta_qT_{(\Lambda^\gamma f)}g\|_{L^2}^2
\lesssim\sum_q\|S_{q-1}\Lambda^\gamma f\cdot\Delta_qg\|_{L^2}^2
\lesssim\sum_q\|S_{q-1}\Lambda^\gamma f\|_{L^2}^2\|\Delta_qg\|_{L^\infty}^2\\
\lesssim&\sum_q2^{2q\gamma}\|S_{q-1}f\|_{L^2}^2\|\Delta_qg\|_{L^\infty}^2
\leq\|f\|_{L^2}^2\|g\|_{B^\gamma_{\infty,2}}^2,\\
\|\RNum{1}_3\|_{L^2}^2\leq&
\sum_q\|\Lambda^\gamma(\tilde{\Delta}_qf\cdot\Delta_qg)\|_{L^2}^2
\lesssim \sum_q2^{2q\gamma}\|\tilde{\Delta}_qf\cdot\Delta_qg\|_{L^2}^2
\leq\sum_q2^{2q\gamma}\|\tilde{\Delta}_qf\|_{L^2}^2\|\Delta_qg\|_{L^\infty}^2
\\
&\leq\|f\|_{L^2}^2\|g\|_{B^\gamma_{\infty,2}}^2,\\
\|\RNum{2}_3\|_{L^2}^2\leq&
\sum_q\|\tilde{\Delta}_q(\Lambda^\gamma f)\cdot\Delta_qg\|_{L^2}^2
\leq \sum_q\|\tilde{\Delta}_q(\Lambda^\gamma f)\|_{L^2}^2\|\Delta_qg\|_{L^\infty}^2
\lesssim\sum_q2^{2q\gamma}\|\tilde{\Delta}_qf\|_{L^2}^2\|\Delta_qg\|_{L^\infty}^2\\
&
\leq \|f\|_{L^2}^2\|g\|_{B^\gamma_{\infty,2}}^2,
\end{align*}
as $C^{\gamma+\epsilon}$ is embedded in $B^\gamma_{\infty,2}$. These
terms are nicely controlled.

The commutator structure is mainly used to estimate
$\RNum{1}_2-\RNum{2}_2$. Let us denote the difference as~\RNum{3}.
Given any $q\in\mathbb{N}$,
\[\Delta_q\RNum{3}=\sum_p\Delta_q\left(\Lambda^\gamma(S_{p-1}g\cdot\Delta_pf)
-S_{p-1}g\cdot\Lambda^\gamma(\Delta_pf)\right)=:\sum_p\RNum{3}_p.\]
Note that $\RNum{3}_p\equiv0$ for $|p-q|\geq5$. Therefore, it is a finite
sum. We discuss $\RNum{3}_q$ and the other terms can be treated similarly.

Followed from \cite{kenig1993well}, we estimate $\RNum{3}_q$ in the Fourier side,
\[
\RNum{3}_q(x)=\iint
(|\xi+\zeta|^\gamma-|\xi|^\gamma)\eta(2^{-q}(\xi+\zeta))\chi(2^{-(q-2)}\zeta)
\eta(2^{-q}\xi)\hat{f}(\xi)\hat{g}(\zeta)e^{i(\xi+\zeta)x} d\xi d\zeta.
\]
Define a multiplier $m(\xi,\zeta)$ as
\[m(\xi,\zeta)=\frac{|\xi+\zeta|^\gamma-|\xi|^\gamma}{|\zeta|^\gamma}
\eta(\xi+\zeta)\chi(4\zeta)\eta(\xi).\]
It is easy to check that $m$ is uniformly bounded, compactly supported and $C^\infty$.
Let $m_q(\xi,\zeta)=m(2^{-q}\xi,2^{-q}\zeta)$, then
\[\RNum{3}_q(x)=\iint
m_q(\xi,\zeta)\hat{f}(\xi)|\zeta|^\gamma\hat{g}(\zeta)e^{i(\xi+\zeta)x}d\xi d\zeta
=\iint h_q(y,z)\cdot\Delta_qf(x-y)\cdot\Lambda^\gamma S_{q-1}g(x-z)dydz,\]
where
\[h_q(y,z)=C\iint
m_q(\xi,\zeta)e^{i(\xi y+\zeta z)}d\xi d\zeta.\]
Compute
\[\iint |h_q(y,z)| dydz=2^{2q}\iint |h_1(2^qy, 2^qz)|dydz=\iint
  |h_1(y,z)|dydz\leq C,\]
where the last integral is bounded due to smoothness of $m$, and the
constant $C$ does not depend on~$q$.
Then, applying Young's inequality, we get
\[\|\RNum{3}_q\|_{L^2}\lesssim\|h_q(\cdot,\cdot)\|_{L^1}\|\Delta_qf\|_{L^2}
\|\Lambda^\gamma S_{q-1}g\|_{L^\infty}\lesssim
\|\Delta_qf\|_{L^2} \sum_{p<q-1}2^{p\gamma}\|\Delta_pg\|_{L^\infty}.\]
We collect all modes and conclude
\begin{align*}
\|\RNum{3}\|_{L^2}^2=&\sum_q\|\Delta_q\RNum{3}\|_{L^2}
\lesssim\sum_q \|\Delta_qf\|_{L^2}^2
\left(\sum_{p< q-1}2^{p\gamma}\|\Delta_pg\|_{L^\infty}\right)^2
\\
\lesssim &\sum_q \|\Delta_qf\|_{L^2}^2\sum_{p<q-1}
2^{2p(\gamma+\frac{\alpha}{2})} \|\Delta_pg\|_{L^\infty}^2\\
=&\sum_p 2^{2p(\gamma+\frac{\epsilon}{2})} \|\Delta_pg\|_{L^\infty}^2
\sum_{q>p+1}\|\Delta_qf\|_{L^2}^2
\leq \|f\|_{L^2}^2\|g\|_{B^{\gamma+\frac{\epsilon}{2}}_{\infty,2}}^2
\lesssim\|f\|_{L^2}^2\|g\|_{C^{\gamma+\epsilon}}^2.
\end{align*}
\end{proof}

%\bibliographystyle{plain}
%\bibliography{fractalBurgers-short}

\end{document}